%% file: cell_vertex_weno_arxiv.tex
\documentclass[preprint, 12pt]{elsarticle}

\input{shared-hdr}

\input{article-hdr}
\input{commands}

\journal{}

\newcommand{\beq}{\begin{equation}}
\newcommand{\eeq}{\end{equation}}

\begin{document}
	
\begin{frontmatter}
	
	\title{Cell-vertex WENO schemes with shock-capturing quadrature for high-order finite element discretizations of hyperbolic problems}
	% Title 2: Quadrature-driven cell-vertex WENO schemes for high-order finite element discretizations of hyperbolic problems
	
	\author[TU]{Joshua Vedral\corref{cor1}}
	\ead{joshua.vedral@math.tu-dortmund.de}
	\cortext[cor1]{Corresponding author}
	\author[TU]{Dmitri Kuzmin}
	\ead{kuzmin@math.uni-dortmund.de}
	
	\address[TU]{Institute of Applied Mathematics (LS III), TU Dortmund University\\ Vogelpothsweg 87, D-44227 Dortmund, Germany}

	\begin{abstract}
          We propose a new kind of localized shock capturing for continuous (CG) and discontinuous Galerkin (DG) discretizations of hyperbolic conservation laws. The underlying framework of dissipation-based weighted essentially nonoscillatory (WENO) stabilization for high-order CG and DG approximations was introduced in our previous work. In this general framework, Hermite WENO (HWENO) reconstructions are used to calculate local smoothness sensors that determine the appropriate amount of artificial viscosity for each cell. In the original version, candidate polynomials for WENO averaging are constructed using the derivative data from von Neumann neighbors. We upgrade this standard `cell-cell' reconstruction procedure by using WENO polynomials associated with mesh vertices as candidate polynomials for cell-based WENO averaging. The Hermite data of individual cells is sent to vertices of those cells, after which vertex-averaged HWENO data is sent back to cells containing the vertices. The new `cell-vertex' averaging procedure includes the data of vertex neighbors without explicitly adding them to the reconstruction stencils. It mitigates mesh imprinting and can also be used in classical HWENO limiters for DG methods. The second main novelty of the proposed approach is a quadrature-driven distribution of artificial viscosity within high-order finite elements. Replacing the linear quadrature weights by their nonlinear WENO-type counterparts, we concentrate shock-capturing dissipation near discontinuities while minimizing it in smooth portions of troubled cells. This redistribution of WENO stabilization preserves the total dissipation rate within each cell and improves local shock resolution without relying on subcell decomposition techniques. Numerical experiments in one and two dimensions demonstrate substantial improvements in accuracy and robustness for high-order elements, providing a compelling alternative to standard cell-cell reconstructions and subcell shock-capturing schemes.
	\end{abstract}
	
	
\end{frontmatter}

% ---------------------------------------------------------------------------
\section{Introduction}
% ---------------------------------------------------------------------------

Extensive research efforts have been invested over the years in the development of high-order non-oscillatory schemes for hyperbolic problems. The essentially non-oscillatory (ENO) polynomial selection principle, introduced by Harten and Osher \cite{harten1987}, laid the foundation for stable high-order reconstructions. Shu and Osher \cite{shu1988,shu1989} further popularized ENO schemes by combining them with high-order Runge-Kutta time stepping and formulating reconstructions directly at the flux level. Building on this framework, Liu et al. \cite{liu1994} introduced the weighted ENO (WENO) methodology, which Jiang and Shu \cite{jiang1996} later equipped with classical smoothness indicators and nonlinear weights, establishing the basis for most subsequent developments.

A significant line of research concerns the integration of WENO ideas into discontinuous Galerkin (DG) schemes. Qiu and Shu \cite{qiu2005} were the first to employ WENO reconstructions as limiters for Runge-Kutta DG (RKDG) schemes, leading to a broad family of DG-WENO methods. In DG approximations, derivative information from neighboring elements is naturally available and can complement element averages in the construction of candidate polynomials. Hermite WENO (HWENO) limiters \cite{qiu2004,qiu2005b} exploit this feature, with extensions to unstructured meshes by Zhu and Qiu \cite{zhu2009,zhu2008}. Luo et al. \cite{luo2007} demonstrated that incorporating first-order derivatives in the construction of candidate polynomials yields a simple and efficient WENO reconstruction on multidimensional unstructured grids. Zhu et al. \cite{zhu2017} later proposed a least-squares variant of the simple WENO limiter \cite{zhong2013}.

Parallel efforts have targeted the WENO reconstruction itself. Henrick et al. \cite{henrick2005} showed that the Jiang-Shu formulation may degrade from fifth- to third-order accuracy at critical points and proposed a simple modification of the nonlinear weights using mapping functions to restore fifth-order accuracy. Borges et al. \cite{borges2008} developed higher-order global smoothness indicators via linear combination of classical ones yielding new non-oscillatory weights.

Recent developments emphasize flexibility and adaptivity of WENO schemes. Central WENO (CWENO) schemes \cite{levy1999,levy2000,qiu2002} split the reconstruction stencil into a centered and a WENO component. WENO schemes with adaptive order (WENO-AO) \cite{balsara2020,balsara2016} hybridize between high-order polynomials on a centered stencil and CWENO reconstructions. Unequal-sized WENO (US-WENO) schemes \cite{zhu2016,zhu2018} employ artificial linear weights to avoid negative weights on unstructured grids without losing accuracy. Targeted ENO (TENO) schemes \cite{fu2016,fu2017a} reduce numerical dissipation by using an ENO-like stencil selection that discards candidates crossed by discontinuities. We refer the interested reader to \cite{navas-montilla2024} for a comparison between classical finite volume WENO and TENO schemes for the Euler equations with gravity. 

Beyond modifications of WENO reconstructions, many high-order shock-capturing methods combine a linear high-order discretization with a nonlinear WENO-type correction. Early examples include the compact upwind–ENO scheme of Adams and Shariff \cite{adams1996} and its WENO counterpart by Pirozzoli \cite{pirozzoli2002}. Hybrid centered-WENO approaches were later proposed in \cite{hill2004,li2020,wan2012}, and central-upwind schemes with WENO shock-capturing for systems of conservation laws were developed in \cite{hu2010}. In fact, the underlying idea of blending linear high-order stabilization with a nonlinear shock-capturing operator predates WENO and can be traced back to the work of Hughes et al. \cite{hughes1987,hughes1986}, Johnson et al. \cite{johnson1990} and Codina \cite{codina1993}. This general stabilization philosophy has since been expanded through nonlinear diffusion and artificial viscosity frameworks, including local projection stabilization with crosswind diffusion \cite{barrenechea2013}, monotonicity-preserving and graph-Laplacian stabilization methods \cite{badia2017a,badia2014}, and residual-based artificial viscosity methods \cite{nazarov2013}. 

Existing HWENO/RKDG approaches associate candidate polynomials almost exclusively with von Neumann (edge/face) neighbors. As a promising alternative, we introduce an averaging procedure in which one candidate polynomial represents the finite element approximation inside a cell, while the remaining ones represent HWENO polynomial reconstructions associated with vertices of the cell. The cell-vertex blending strategy minimizes mesh imprinting effects by incorporating vertex-neighbor information into the final cell averages. The new approach can be used both in the context of RKDG schemes equipped with HWENO limiters and in the framework of dissipation-based WENO stabilization via nonlinear shock-capturing viscosity. The latter approach was introduced in \cite{kuzmin2023a} in the continuous Galerkin (CG) setting and extended to DG schemes in \cite{vedral-arxiv}. As a further improvement, we propose a novel quadrature-based redistribution of artificial dissipation within each cell. The weights of a standard quadrature rule are adjusted similarly to the linear weights of a WENO reconstruction, and artificial viscosity is redistributed accordingly. This enhancement of our dissipation-based WENO finite element schemes has the same effect as a localization of nonlinear stabilization to subcells \cite{glaubitz2019,hennemann2021,markert2023,rueda-ramirez2022,vilar2024}, but is simpler and more general. 

We continue in Section \ref{sec:stab} with a review of the dissipation-based WENO stabilization methodology. The new cell-vertex reconstruction procedure and the proposed quadrature-driven localization of shock-capturing viscosity are introduced in Sections \ref{sec:cell-vertex} and \ref{sec:sc-quad}, respectively. In Section \ref{sec:results}, we provide further implementation details and run numerical experiments in one and two space dimensions. The main findings are summarized in Section \ref{sec:concl}.

% ---------------------------------------------------------------------------
\section{Stabilized finite element formulation}
\label{sec:stab}
% ---------------------------------------------------------------------------

Let $u(\mathbf{x},t)$ be a conserved quantity or a vector of conserved variables depending on the space location $\mathbf{x}$ and time instant $t\ge 0$. For simplicity, we impose periodic boundary conditions on the boundary $\Gamma$ of a spatial domain $\Omega\subset \mathbb{R}^d$, $d\in\{1,2,3\}$. We consider the initial value problem
\begin{subequations}
	\begin{alignat}{2}
		\frac{\partial u}{\partial t}+\nabla\cdot\mathbf{f}(u)&=0 \quad &&\text{in }\Omega \times (0,T], \label{eq:pde}\\
		u(\cdot,0)&=u_0 \quad &&\text{in } \Omega,
	\end{alignat}\label{eq:ivp}
\end{subequations}
where $\mathbf{f}(u)$ is the flux function, $T>0$ is the final time, and $u_0$ is the initial datum.

\subsection{Galerkin scheme}

We discretize \eqref{eq:ivp} in space using a CG or DG method on a conforming, possibly unstructured, affine mesh $\mathcal{T}_h$ that consists of $d$-dimensional simplices or boxes $K_e$, $e=1,\ldots,E_h$ such that $\bar{\Omega}=\cup_{e=1}^{E_h}K_e$. We denote the mesh size by $h=\max_{K_e\in\mathcal{T}_h}h_e$, where $h_e=\diam(K_e)$. The corresponding finite element spaces for DG and CG discretizations are defined as
\begin{align*}
	V_{h,p}^d(\mathcal{T}_h)&:=\{v_h\in L^2(\Omega):v_h|_{K_e}\in V_{h,p}(K_e)\,\forall K_e \in \mathcal{T}_h\}, \\ V_{h,p}^c(\mathcal{T}_h)&:=V_{h,p}^d(\mathcal{T}_h)\cap C^0(\bar{\Omega}),
\end{align*}
respectively. Here, $V_{h,p}(K_e)\in\{\mathbb{P}_p(K_e),\mathbb{Q}_p(K_e)\}$ denotes the local polynomial space of (total) degree up to $p\in \mathbb{N}$. For better readability, we omit the index $p$, as well as the superscripts $d$ or $c$, in what follows. The simplified notation $V_h$ will refer to either $V_{h,p}^d$ or $V_{h,p}^c$, depending on the type of the Galerkin discretization. We then seek an approximate solution of the form
\begin{equation}
	u_h(\mathbf{x},t)=\sum_{j=1}^{N_h}u_j(t)\varphi_j(\mathbf{x}),
\end{equation}
where $(\varphi_j)_{j=1}^{N_h}$ is any nodal basis for $V_h$ with polynomial degree $p$. The associated nodal points are denoted by $\mathbf{x}_1,\ldots,\mathbf{x}_{N_h}$.

The standard Galerkin discretization of \eqref{eq:ivp} yields the semi-discrete problem
\begin{equation}
	\sum_{e=1}^{E_h}\int_{K_e}v_h\Bigg(\frac{\partial u_h}{\partial t}+\nabla \cdot \mathbf{f}(u_h)\Bigg)\,\mathrm{d}\mathbf{x}=0\quad \forall v_h\in V_h.
	\label{eq:galsd1}
\end{equation}
Applying integration by parts, we rewrite \eqref{eq:galsd1} in the generic form
\begin{equation}
	\frac{\mathrm{d}}{\mathrm{dt}}(u_h,v_h)+a(u_h,v_h)=0 \quad \forall v_h\in V_h,
	\label{eq:galsd2}
\end{equation}
where
\begin{equation}
	(u_h,v_h)=\sum_{e=1}^{E_h}\int_{K_e}v_hu_h\,\mathrm{d}\mathbf{x}
\end{equation}
and
\begin{equation}
	a(u_h,v_h) = -\sum_{e=1}^{E_h}\int_{K_e}\nabla v_h \cdot \mathbf{f}(u_h)\,\mathrm{d}\mathbf{x}+\sum_{e=1}^{E_h}\int_{\partial K_e}v_h\mathbf{f}(u_h)\cdot\mathbf{n}\,\mathrm{ds}.
	\label{eq:convop}
\end{equation}
For CG methods, the second term in \eqref{eq:convop} vanishes due to continuity across element interfaces. For DG methods, the solution may be discontinuous at element boundaries, and the flux $\mathbf{f}(u_h)\cdot \mathbf{n}$ must be replaced by a Riemann solver. In this work, we use the local Lax-Friedrichs flux for scalar problems and the Harten-Lax-van Leer (HLL) flux for hyperbolic systems. Estimates of the maximum wave speed for the Euler equations can be found in \cite{guermond2016a}. 

\subsection{Dissipative nonlinear stabilization}

To ensure stability and achieve optimal convergence rates, we equip the standard Galerkin scheme \eqref{eq:galsd2} with the stabilization operator of the dissipation-based WENO framework \cite{kuzmin2023a,vedral-arxiv}. The stabilized formulation reads
\begin{equation}
	\frac{\mathrm{d}}{\mathrm{dt}}(u_h,v_h)+a(u_h,v_h)+s_h(u_h,v_h)=0 \quad \forall v_h\in V_h,
	\label{eq:stabform}
\end{equation}
where the operator $s_h$ is assembled from local element contributions as follows:
\begin{equation}
	s_h(u_h,v_h)=\sum_{e=1}^{E_h}s_h^e(u_h,v_h).
\end{equation}
We emphasize that this construction leaves the numerical fluxes unchanged.

The stabilization must (i) suppress spurious oscillations near discontinuities and (ii) preserve optimal accuracy in smooth regions. In the classical JST scheme \cite{jameson1981}, this is achieved by blending high- and low-order stabilization terms. The finite element version proposed in \cite{kuzmin2023a} uses
\begin{equation}
	s_h^e(u_h,v_h)=\gamma_es_h^{e,H}(u_h,v_h)+(1-\gamma_e)s_h^{e,L}(u_h,v_h).
	\label{eq:stab}
\end{equation}
Here, the parameter $\gamma_e\in [0,1]$ serves as an elementwise smoothness indicator.

Following \cite{kuzmin2023a,vedral-arxiv}, we define the low-order operator as
\begin{equation}
	s_h^{e,L}(u_h,v_h)=\nu_e\int_{K_e}\nabla v_h\cdot \nabla u_h \, \mathrm{d}\mathbf{x},
	\label{eq:lostab}
\end{equation}
where $\nu_e$ is the viscosity parameter given by $\nu_e=\frac{\lambda_eh_e}{2p}$ and $\lambda_e=\|\mathbf{f}(u_h)\|_{L^{\infty}(K_e)}$ denotes the local maximum wave speed. 

The construction of the high-order stabilization operator is somewhat more delicate.
Adding an antidiffusive correction to the low-order stabilization operator \eqref{eq:lostab}, we obtain its high-order counterpart
\begin{equation}
	s_h^{e,H}(u_h,v_h)=\nu_e\int_{K_e}\nabla v_h\cdot(\nabla u_h - P_h\nabla u_h)\,\mathrm{d}\mathbf{x},
	\label{eq:hostab}
\end{equation}
where $P_h:L^2(\Omega)\to V_h$ is a projection operator. As shown in \cite{vedral2025}, the choice of $P_h$ is critical for achieving optimal convergence rates.

For CG discretizations, we choose $P_h$ to be the (componentwise) Scott-Zhang variant \cite{scott1990} of the Cl{\'e}ment quasi-interpolation operator \cite{clement1975}. For DG methods, we use the (componentwise) $L^2$ orthogonal projection. In this case, $P_h\nabla u_h = \nabla u_h$ on each element, and hence $s_h^{e,H}=0$. No high-order stabilization is required for DG methods equipped with stable numerical fluxes. For CG methods, however, some form of high-order stabilization is essential to obtain optimal convergence rates; see, e.g., \cite{quarteroni1994}.
\begin{remark}
	The operator \eqref{eq:hostab} is closely related to several existing stabilization techniques. In particular, it can be viewed as a component of a two-level variational multiscale (VMS) method \cite{john2006,kuzmin2020f,lohmann2017} and as a local fluctuation operator acting on the solution gradient, similar to those used in local projection stabilization (LPS) methods \cite{braack2006,barrenechea2010,knobloch2010a,knobloch2009}. It is also closely related to the orthogonal subscale VMS approach of Codina and Blasco \cite{codina1997,codina2002}.
\end{remark}
\begin{remark}
	For CG discretizations, the blended operator \eqref{eq:stab} can be written in the simplified form \cite{kuzmin2025a}
	\begin{equation}
		s_h^e(u_h,v_h)=\nu_e\int_{K_e}\nabla v_h\cdot (\nabla u_h - \gamma_e P_h\nabla u_h)\,\mathrm{d}\mathbf{x}.
	\end{equation}
\end{remark}

The blending factor $\gamma_e$ controls the interplay between the high- and low-order stabilization terms. When $\gamma_e=0$, the nonlinear form $s_h^e$ in \eqref{eq:stab} reduces to a low-order local Lax-Friedrichs-type diffusion. In contrast, when $\gamma_e=1$, it recovers the high-order operator \eqref{eq:hostab}. For linear finite elements and the lumped-mass $L^2$ projection operator $P_h$, the bilinear form $s_h^{e,H}(u_h,v_h)$ is symmetric and positive semi-definite \cite{olshanskii2025}. The consistent-mass $L^2$ projection provides these properties for any polynomial degree $p\ge 1$.

The smoothness sensor proposed in \cite{kuzmin2023a} is given by
\begin{equation}
	\gamma_e = 1-\min\Bigg(1,\frac{\|u_h^e-u_h^{e,\ast}\|_e}{\|u_h^e\|_e}\Bigg)^q.
	\label{eq:gamma}
\end{equation}
We adopt this definition of $\gamma_e$, which
quantifies deviations of $u_h^e:=u_h|_{K_e}$ from a WENO reconstruction $u_h^{e,\ast}$. The exponent $q$ controls the sensitivity of $\gamma_e\in [0,1]$ to the relative discrepancy in the partial derivatives. The deviations are measured using the scaled Sobolev semi-norm \cite{friedrich1998,jiang1996}
\begin{equation}
	\|v\|_e=\Bigg(\sum_{1\le |\mathbf{k}|\le p}h_e^{2|\mathbf{k}|-d}\int_{K_e}|D^{\mathbf{k}}v|^2\,\mathrm{d}\mathbf{x}\Bigg)^{1/2} \quad \forall v\in H^p(K_e),
	\label{eq:snorm}
\end{equation}
where $\mathbf{k}=(k_1,\ldots,k_d)$ denotes the multiindex of the partial derivative $D^{\mathbf{k}}v$.

By replacing the test function $v_h$ with the basis function $\varphi_i$ in the stabilized formulation \eqref{eq:stabform}, we obtain the semi-discrete system 
\begin{equation}
  	\frac{\mathrm{d}}{\mathrm{dt}}
	(u_h,\varphi_i)+a(u_h,\varphi_i)+s_h(u_h,\varphi_i)=0,\qquad
        i=1,\ldots,N_h.
	\label{eq:sdstab}
\end{equation}
The numerical solution can then be advanced in time using, e.g., a strong stability preserving (SSP) Runge-Kutta method \cite{gottlieb2001}.

\begin{remark}
An element-based monolithic convex limiting (MCL) scheme for the CG version of \eqref{eq:sdstab}
was recently developed in \cite{kuzmin2025a} to ensure positivity
preservation. The DG version can be limited in the same way.
\end{remark}

% ---------------------------------------------------------------------------
\section{Cell-vertex WENO reconstruction}
\label{sec:cell-vertex}
% ---------------------------------------------------------------------------

In principle, any standard WENO scheme from the literature can be used to reconstruct the
polynomial $u_h^{e,\ast}$ for smoothness estimation using the shock detector \eqref{eq:gamma}.
In our previous work \cite{kuzmin2023a,vedral-arxiv,vedral2025}, we adopted the Hermite WENO (HWENO) framework of Qiu and Shu \cite{qiu2004} and introduced several modifications. Unlike HWENO formulations that use only cell averages and first-order derivatives \cite{luo2007,qiu2004}, our approach incorporates all partial derivatives up to order $p$, resulting in inherently compact stencils.

The procedure proposed in this section is based on the same design principles, but the candidate polynomials corresponding to finite element approximations in edge/face neighbors of a cell are replaced by WENO polynomials associated with vertices of that cell. These auxiliary polynomials are constructed using the data from cells containing the vertex. The new reconstruction procedure can be used to construct $u_h^{e,\ast}$ for \eqref{eq:gamma} or for an HWENO limiter that overwrites the DG approximation $u_h^{e}$ by $u_h^{e,\ast}$.

\subsection{Cell-cell WENO reconstruction}

One of the key challenges in designing accurate and robust WENO schemes is the choice of reconstruction stencils. In finite volume (FV) methods, high-order polynomials are fitted to the cell averages of approximate solutions in selected neighbor cells. As noted in \cite{luo2007, tsoutsanis2019}, the number of cells required to construct a polynomial of degree $p_r$ in $d$ space dimensions is
\begin{equation}
	N_c^{\triangle} = \frac{(p_r+1)(p_r+2)\cdots(p_r+d)}{d!}, \quad N_c^{\square} = (p_r+1)^d
\end{equation}
for simplicial and tensor-product cells, respectively. Choosing appropriate stencils becomes increasingly difficult for high-order polynomials and in higher dimensions. The task becomes even more challenging on unstructured grids, which impose additional constraints on the selection of stencil cells. A detailed study of central stencil algorithms in the FV setting can be found in \cite{tsoutsanis2019}.

In the DG framework, standard WENO reconstructions from element averages of DG polynomials require FV-type stencils, which leaves the problem of identifying suitable cells unresolved. A key advantage of DG (and CG) methods is that derivative information is readily available and can be exploited to reduce the stencil size. This idea leads to the class of Hermite WENO (HWENO) schemes, as discussed in \cite{luo2007,qiu2004,qiu2005b,zhu2009}.

The $l$th candidate polynomial $u_{h,l}^e$ of an HWENO reconstruction
for a target cell $K_e$ uses the data of a stencil cell $K_{e'}$. The conditions
\begin{equation}
	\frac{1}{|K_e|}\int_{K_e}u_{h,l}^e\,\mathrm{d}\mathbf{x}=\frac{1}{|K_e|}\int_{K_e}u_h^e\,\mathrm{d}\mathbf{x}, \quad D^{\mathbf{k}}u_{h,l}^e=D^{\mathbf{k}}u_h^{e'}, \quad 1\le |\mathbf{k}|\le p_r
	\label{eq:hweno1}
\end{equation}
are satisfied by the Hermite interpolation polynomials
\begin{equation}
	u_{h,l}^e(\mathbf{x})=u_h^{e'}(\mathbf{x})+\frac{1}{|K_e|}\int_{K_e}(u_h^e-u_h^{e'})\,\mathrm{d}\mathbf{x}, \quad \mathbf{x}\in K_e.
	\label{eq:hweno2}
\end{equation}
In this \emph{cell-cell} version of the reconstruction procedure,
each candidate polynomial is associated with a cell $K_{e'}$
that shares an edge/face $\partial K_e\cap \partial K_{e'}$
with $K_e$. Such cells constitute the von Neumann neighborhood of $K_e$.

\subsection{Cell-vertex WENO scheme}
\label{subsec:cvweno}

The classical cell-cell construction of candidate polynomials in \eqref{eq:hweno1} and \eqref{eq:hweno2} has a key limitation: it only uses information from von Neumann neighbors and ignores vertex neighbors. The reconstruction stencils are directionally biased, because information propagates only across faces. This can be particularly problematic for oblique shocks or features not aligned with the grid. As a result, smooth regions may be misidentified, and the overall accuracy may deteriorate. A straightforward fix would be to directly include all vertex neighbors in the reconstruction stencil, which, however, increases the computational cost.

To address the above issues, we propose a vertex-based WENO reconstruction procedure, referred to as a \textit{cell-vertex} WENO scheme. This approach incorporates information from vertex neighbors while keeping the number of candidate polynomials similar to that for traditional cell-cell approaches. The reconstruction on each element $K_e$ is carried out in two stages. First, vertex-based candidate polynomials are constructed using WENO averaging of cell data. Second, the WENO reconstruction $u_h^{e,\ast}$ is defined by blending $u_h^{e}$ with vertex candidates. The corresponding weights are calculated following the standard Jiang-Shu methodology. In the remainder of this section, we describe the details of our cell-vertex reconstruction procedure.

Let $K_e$ be a mesh cell with vertices $\mathbf{x}_j\in \mathcal{V}(K_e)$, where $\mathcal{V}(K_e)$ denotes the set of all vertices of $K_e$, and let $v_e:=|\mathcal{V}(K_e)|$ be the number of cell vertices. For each vertex $\mathbf{x}_i$, we define its vertex patch $\mathcal{E}(\mathbf{x}_i):=\{K_{e'}:\mathbf{x}_i\in K_{e'}\}$ and store the indices of elements that share $\mathbf{x}_i$
in $\mathcal C(\mathbf x_i):=\{e'\,:\,K_{e'}\in\mathcal{E}(\mathbf{x}_i)\}$.

In our cell-vertex scheme, we combine the cell candidate polynomial
\begin{equation}
	u_{h,0}^e:=u_h^e
\end{equation} 
with $v_e$ candidate polynomials  $u_{h,l}^e=u_{h,i}^*$ that are associated with the vertices $\mathbf{x}_l^e=\mathbf x_i,\ i\in\{1,\ldots,N_h\}$ of $K_e$. The same vertex-based WENO reconstruction $u_{h,i}^*$ serves as a candidate polynomial for all $K_e$, $e\in \mathcal{C}(\mathbf{x}_i)$. The vertex candidates $u_{h,l}^e$ of a cell $K_e$ are numbered using the local index $l\in\{1,\ldots,v_e\}$.

\begin{figure}[t]
	\centering
	\includegraphics[width=0.23\textwidth]{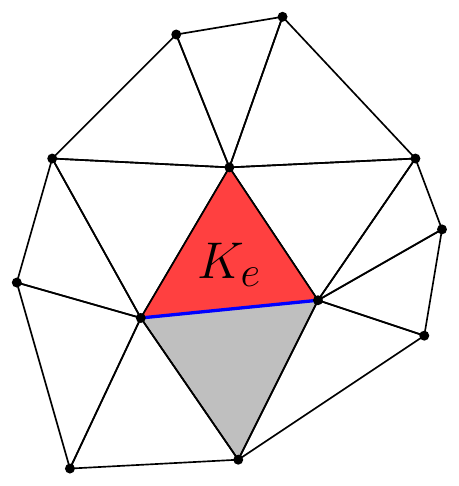}
	\hfill
	\includegraphics[width=0.23\textwidth]{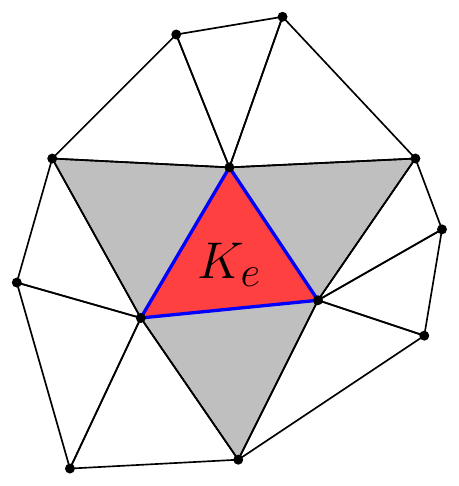}
	\hfill
	\includegraphics[width=0.23\textwidth]{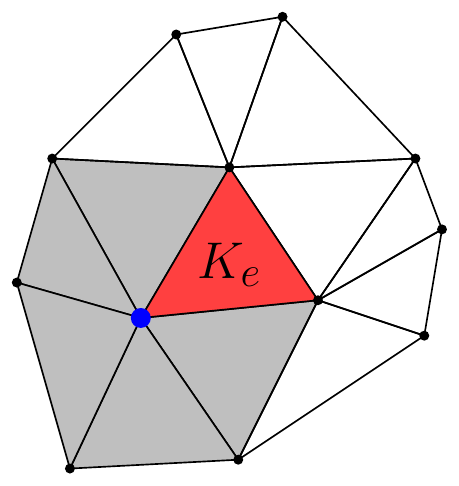}
	\hfill
	\includegraphics[width=0.23\textwidth]{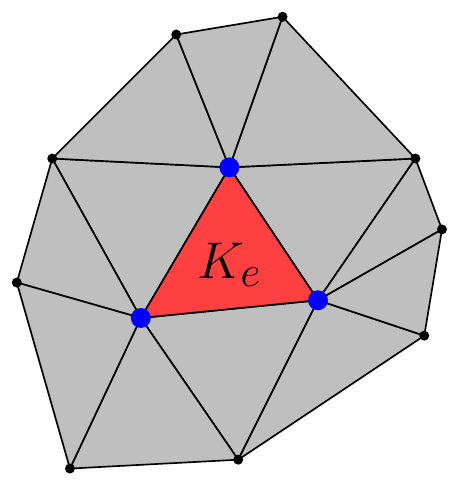}
	\caption{Candidate polynomial stencils for triangular elements. The left two figures show cell-cell stencils. The first one illustrates the cell-based stencil for a single candidate polynomial, while the second one displays the full stencil of a cell-cell reconstruction. The right two figures present the corresponding vertex-based and cell-vertex stencils. In both approaches, four candidate polynomials are used, including the finite element approximation $u_h^e$ in the target cell $K_e$.}
	\label{fig:stencils}
\end{figure}

To construct the vertex-based polynomials $u_{h,i}^*$ in a loop over
$i=1,\ldots,N_h$, we exploit the fact that the restriction of the finite element solution to each element $K_{e'}\in\mathcal{E}(\mathbf{x}_i)$ coincides with the Taylor expansion \cite{kuzmin2010,luo2008}
\begin{equation}
	u_h^{e'}(\mathbf{x})=\sum_{0\le |\mathbf{k}|\le p}\frac{\partial^{\mathbf{k}}u_h^{e'}(\mathbf{x}_i)}{\mathbf{k}!}(\mathbf{x}-\mathbf{x}_i)^{\mathbf{k}}.
	\label{eq:taylorsol}
\end{equation}
The vertex-based candidate polynomial $u_{h,l}^e$ for
$\mathbf{x}_l^e=\mathbf x_i$ is then constructed by
taking a weighted sum of these expansions over the vertex patch: 
\begin{equation}
	u_{h,i}^*(\mathbf{x})= \sum_{e'\in\mathcal C(\mathbf x_i)}\omega_{e'}^iu_h^{e'}(\mathbf{x})=\sum_{0\le |\mathbf{k}|\le p}\frac{(\mathbf{x}-\mathbf{x}_i)^{\mathbf{k}}}{\mathbf{k}!}\sum_{e'\in\mathcal C(\mathbf x_i)}\omega_{e'}^i\partial^{\mathbf{k}}u_h^{e'}(\mathbf{x}_i).
	\label{eq:taylorcand}
\end{equation}
Thus, $u_{h,i}^*$ is a convex combination of Taylor polynomials representing the
finite element solution in cells around $\mathbf{x}_i$. The stencils used in the cell-based and vertex-based polynomial reconstructions are illustrated in Fig.~\ref{fig:stencils}.

The candidate polynomials defined in \eqref{eq:taylorcand} do not necessarily preserve the cell average of the finite element solution on $K_e$. This is not an issue, since the smoothness indicator $\gamma_e$ in \eqref{eq:gamma} depends only on the derivatives of the candidate polynomials. However, the cell averages can be corrected similarly to \eqref{eq:hweno2}. The resulting candidate polynomials, which now preserve the cell average on $K_e$, read
\begin{equation}
	\bar{u}_{h,i}^*(\mathbf{x})= \sum_{e'\in\mathcal C(\mathbf x_i)}\omega_{e'}^iu_h^{e'}(\mathbf{x})+\frac{1}{|K_e|}\int_{K_e}(u_h^e-\sum_{e'\in\mathcal C(\mathbf x_i)}\omega_{e'}^iu_h^{e'})\,\mathrm{d}\mathbf{x}.
\end{equation}

The nonlinear weights for vertex-based WENO averaging in \eqref{eq:taylorcand} are calculated using equal linear weights in the Jiang-Shu formula \cite{jiang1996} such that
\begin{equation}
	\omega_{e'}^i=\frac{\tilde{\omega}_{e'}^i}{\sum_{e\in\mathcal C(\mathbf x_i)}\tilde{\omega}_e^i}, \quad \tilde{\omega}_{e}^i=\frac{1}{(\varepsilon+\beta_{e}^i)^r},
\end{equation}
where $\varepsilon>0$ prevents division by zero, $\beta_{e}^i$ are pointwise smoothness indicators, and $r>0$ is a positive integer. Following Jiang and Shu \cite{jiang1996}, we set $\varepsilon=10^{-6}$ and $r=2$. The smoothness indicators are given by 
\begin{equation}
	\beta_{e}^i=\|u_h^{e}\|_{\mathbf x_i}^q, \quad q\ge 1,
\end{equation}
where the vertex-based semi-norm
\begin{equation}
	\|v\|_{\mathbf x_i}=\Bigg(\sum_{1\le |\mathbf{k}|\le p}h_e^{2|\mathbf{k}|}|D^{\mathbf{k}}v(\mathbf{x}_i)|^2\Bigg)^{1/2} \quad \forall v\in H^s(K_e)
	\label{eq:pointsnorm}
\end{equation}
is the pointwise analogue of the element-level semi-norm in \eqref{eq:snorm}. By the Sobolev embedding theorem, this semi-norm is well-defined for $s>p+d/2$.

Having constructed all vertex-based candidate polynomials $u_{h,i}^*$ in the loop
over $i=1,\ldots,N_h$, we enter a loop over the cell indices $e=1,\ldots,E_h$. In
this loop, we calculate the final cell-vertex WENO reconstructions
\begin{equation}
	u_h^{e,\ast}=\sum_{l=0}^{v_e}\omega_l^eu_{h,l}^e\in \mathbb{P}_p(K_e).
	\label{eq:finalweno}
\end{equation}
The cell candidate $u_{h,0}^e=u_h^e$ and vertex candidates $u_{h,l}^e,\
l=1,\ldots,v_e$ corresponding to $u_{h,i}^*,\ i=i(l)\in\mathcal V(K_e)$
are averaged using the nonlinear weights
\begin{equation}
	\omega_l^e=\frac{\tilde{\omega}_l^e}{\sum_{j=0}^{v_e}\tilde{\omega}_j^e}, \quad \tilde{\omega}_l^e=\frac{\gamma_l^e}{(\varepsilon+\beta_l^e)^r}.
\end{equation}
In the above definition of $\tilde{\omega}_l^e$, the smoothness indicators \cite{friedrich1998,jiang1996}
\begin{equation}
	\beta_l^e=\|u_{h,l}^e\|_e^q, \quad q\ge1
\end{equation}
use the element-level semi-norm from \eqref{eq:snorm}. Following standard practice, we again set $\varepsilon=10^{-6}$ and $r=2$. The linear weights $\gamma_l^e$ can be any positive numbers that sum to one. We assign small linear weights $\gamma_l^e=10^{-3}$ to the vertex-based candidates $u_{h,l}^e$ for $l=1,\ldots,v_e$, while the cell candidate $u_{h,0}^e=u_h^e$ receives the dominant weight $\gamma_0^e=1-\sum_{l=1}^{v_e}\gamma_l^e$, ensuring that the reconstruction reduces to the finite element solution in smooth regions.

\begin{remark} % More detailed?
  When it comes to evaluating the WENO reconstruction \eqref{eq:finalweno} at $\mathbf x\in K_e$, we transform the Taylor basis representation \eqref{eq:taylorcand} of each vertex candidate $u_{h,l}^e=u_{h,i}^*$ into the finite element basis used for the representation of the cell candidate $u_{h,0}^e$. Using the same basis for all $v_e+1$ candidates is particularly convenient for computing the smoothness sensor $\gamma_e$ in \eqref{eq:gamma}.
\end{remark}

\begin{remark}
  The size of candidate stencils that we use in each stage of the cell-vertex approach is comparable to that of a single-stage cell-cell reconstruction. Specifically, both versions use three candidates in one dimension, four candidates for triangles, and five candidates for quadrilaterals and tetrahedra. For cubes, the vertex-based approach slightly increases the number of candidates from six to eight. On unstructured meshes, a minimum-angle condition imposes an upper bound on the number of elements that may meet at a vertex. In contrast to cell-cell reconstructions  \cite{tsoutsanis2019,zou2025}, the cell-vertex stencils  are uniquely determined by the connectivity pattern of the mesh, and the WENO averaging can be coded using standard data structures.
\end{remark}

% ---------------------------------------------------------------------------
\section{Shock-capturing WENO quadrature}
\label{sec:sc-quad}
% ---------------------------------------------------------------------------

In the dissipation-based WENO framework of Section ~\ref{sec:stab}, the smoothness sensor $\gamma_e$ controls the amount of artificial viscosity added to the Galerkin discretization. Since $\gamma_e$ is constant on each element, the corresponding dissipation is distributed uniformly on $K_e$. For high-order elements on coarse meshes, however, it is often desirable to localize the dissipation within an element, concentrating it near sharp gradients or discontinuities. 
We propose a novel WENO-inspired strategy that accomplishes this by rescaling the weights of the quadrature rule for calculating the stabilization term \eqref{eq:stab}. Unlike subcell-based localization approaches \cite{hajduk2020,kuzmin2020a,persson2006,sonntag2014}, our method requires no additional subcell decompositions/computations. It simply redistributes the total dissipation added in $K_e$ among the quadrature points.

Let $D_e=s_h^e(u_h,u_h)$ denote the elementwise dissipation rate for $s_h^e$ defined by \eqref{eq:stab}. In the DG setting, $D_e$ can be written as
\begin{equation}
	D_e^{\text{DG}}=(1-\gamma_e)\nu_e\int_{K_e}|\nabla u_h|^2\,\mathrm{d}\mathbf{x}=\nu_e\int_{K_e}\big(|\nabla u_h|^2-\gamma_e|\nabla u_h|^2\big)\,\mathrm{d}\mathbf{x}
\end{equation}
and in the CG setting as
\begin{equation}
	D_e^{\text{CG}}=\nu_e\int_{K_e}\big(|\nabla u_h|^2-\gamma_e\nabla u_h\cdot P_h\nabla u_h\big)\,\mathrm{d}\mathbf{x}.
\end{equation}

If $s_h^e(u_h,v_h)$ is calculated
using a quadrature rule with weights $\omega_q$ and nodes $\mathbf{x}_q$, $q=1,\ldots,N_q$, both expressions take the unified form
\begin{equation}
	D_e^{\text{orig}}=\nu_e\sum_{q=1}^{N_q}\omega_qf_q,
	\label{eq:dissrateorig}
\end{equation}
where $f_q=(1-\gamma_e)|\nabla u_h(\mathbf{x}_q)|^2$ in the DG setting and
$f_q=|\nabla u_h(\mathbf{x}_q)|^2-\gamma_e\nabla u_h\cdot P_h\nabla u_h$ in the CG setting.

To localize the WENO dissipation, we introduce nonlinear scaling factors $\alpha_q$ and calculate $s_h^e(u_h,v_h)$ using the weights $\alpha_qw_q$ instead of $w_q$. This adjustment of the quadrature rule has the effect of replacing  $D_e^{\text{orig}}$ with
\begin{equation}
	D_e^{\text{scaled}}=\nu_e\sum_{q=1}^{N_q}\omega_q\alpha_qf_q.
	\label{eq:dissratescaled}
\end{equation}
The total dissipation must remain unchanged. This is enforced by defining
\begin{equation}
	\alpha_q:=\frac{\tilde{\alpha}_q\sum_{r=1}^{N_q}\omega_{r}f_r}{\sum_{r=1}^{N_q}\omega_{r}\tilde{\alpha}_{r}f_r},
	\label{eq:alphaq}
\end{equation}
where the preliminary weights
\begin{equation}\label{eq:alpha_prelim}
	\tilde{\alpha}_q:=\frac{\|u_h\|_{\mathbf{x}_q} }{\|u_h\|_e}
\end{equation}
serve as WENO-type smoothness indicators that allocate more dissipation to quadrature nodes near sharp features.

The above choice of the redistribution weights $\alpha_q$ preserves the dissipation rate $D_e^{\text{orig}}$ exactly. Indeed, substituting \eqref{eq:alphaq} into
	\eqref{eq:dissratescaled} gives
	\begin{equation}
		D_e^{\text{scaled}}=\nu_e\sum_{q=1}^{N_q}\omega_q\Bigg(\frac{\tilde{\alpha}_q\sum_{r=1}^{N_q}\omega_rf_r}{\sum_{r=1}^{N_q}\omega_r\tilde{\alpha}_rf_r}\Bigg)f_q=\frac{\nu_e\sum_{r=1}^{N_q}\omega_rf_r}{\sum_{r=1}^{N_q}\omega_r\tilde{\alpha}_rf_r}\sum_{q=1}^{N_q}\omega_q\tilde{\alpha}_qf_q.
	\end{equation}
Using the identity $\nu_e
        \sum_{r=1}^{N_q}\omega_rf_r=D_e^{\text{orig}}$, we find that
	\begin{equation}
	  D_e^{\text{scaled}}=\frac{D_e^{\text{orig}}}{\sum_{r=1}^{N_q}\omega_r\tilde{\alpha}_rf_r}\sum_{q=1}^{N_q}\omega_q\tilde{\alpha}_qf_q=D_e^{\text{orig}}.
	\end{equation}
Furthermore, the adjustment of the quadrature weights does not affect the zero sum property $s_h^e(u_h,1)=0$ that ensures discrete conservation.

The redistribution affects the solution only when a noticeable level of artificial viscosity is present. For smooth elements, where $\gamma_e$ is small and the added dissipation is negligible, the scaling has little impact on the results. In practice, it is therefore worthwhile to  redistribute only when $D_e^{\text{orig}}>C_{D_e^{\text{orig}}}$ or $\gamma_e<C_{\gamma_e}$, where $C_{D_e^{\text{orig}}}$ and $C_{\gamma_e}$ are prescribed cutoff parameters. 

\begin{figure}[t] % Is this figure necessary?
	\centering
	\includegraphics[width=0.96\textwidth]{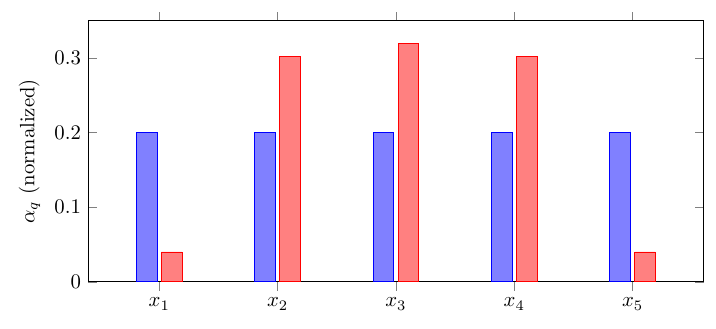}
	\caption{Normalized scaling factors $\alpha_q/\sum_{r=1}^{N_q}\alpha_r$ for a one-dimensional fifth-order finite element containing a discontinuity at its center. The uniform scaling (blue) uses $\alpha_q=1$. The WENO scaling (red) uses the weights $\alpha_q$ given by \eqref{eq:alphaq}.}
	\label{fig:quadweights}
\end{figure}

\begin{remark}
For linear finite elements, the gradient of $u_h$ is piecewise constant. Therefore,
the preliminary weights \eqref{eq:alpha_prelim} are all equal to $\tilde \alpha_q=1/|K_e|$
and the WENO weights \eqref{eq:alphaq} are all equal to unity. Hence, no redistribution is
performed in the case $p=1$. For a finite element of degree $p>1$, the gradient is a
polynomial of degree $p-1$. Hence, the distribution of normalized scaling factors is
nonuniform, as shown in Fig.~\ref{fig:quadweights}. The WENO formula \eqref{eq:alphaq}
assigns larger weights $\alpha_q$ to quadrature points $\mathbf x_q$ at which the
pointwise derivative sensor $\|u_h\|_{\mathbf{x}_q}$ is large compared to the cell average 
 $\|u_h\|_e/|K_e|$.
\end{remark}

% ---------------------------------------------------------------------------
\section{Numerical results}
\label{sec:results}
% ---------------------------------------------------------------------------

In this section, we present numerical results for standard hyperbolic test problems, including linear and nonlinear scalar equations and the Euler equations of gas dynamics. We assess the shock-capturing capabilities of our method as well as its accuracy in smooth solution regimes for both continuous and discontinuous finite element discretizations. Moreover, we verify that the proposed scheme does not produce entropy-violating solutions and highlight its advantages over existing high-order \textit{cell-cell} WENO schemes, particularly for large polynomial degrees.

For clarity, we refer to the \textit{cell-cell} WENO schemes from \cite{kuzmin2023a, vedral2025} as CC-WENO. The \textit{cell-vertex} WENO scheme proposed in Sec.~\ref{sec:cell-vertex} is denoted by CV-WENO. We write CV-WENO-SC for the CV-WENO scheme equipped with shock-capturing WENO quadrature, as introduced in Sec.~\ref{sec:sc-quad}.

All numerical experiments are performed on uniform structured meshes. Extensions to unstructured meshes and three-dimensional problems will be reported elsewhere. The steepening parameter in \eqref{eq:gamma} is set to $q=1$. To avoid unnecessary rescaling of the weights in the WENO quadrature procedure, we only employ the shock-capturing WENO quadrature when $\gamma_e<C_{\gamma_e}$ with $C_{\gamma_e}=0.9$. We use $\mathbb{Q}_p$ Lagrange finite elements with polynomial degree $p \ge 1$. To enable a fair comparison across different polynomial orders, the total number of degrees of freedom $N_h$ is kept fixed by coarsening the mesh as $p$ increases.

Spatial discretization errors are measured using the $L^1$ norm and experimental orders of convergence (EOC) are computed following \cite{leveque1996,lohmann2017}. Time integration is performed using the third-order strong stability preserving (SSP) Runge-Kutta method \cite{gottlieb2001} with time steps chosen sufficiently small so that temporal errors are negligible relative to spatial discretization errors.

All computations are carried out using the C++ finite element library MFEM \cite{mfem2021, mfem2024, mfem} and the results to two-dimensional problems are visualized using the C++ software GLVis \cite{glvis}.

\begin{subsection}{Solid body rotation}
	In our first numerical experiment, we consider LeVeque's solid body rotation problem \cite{leveque1996}, which is used to assess a method's capability to capture smooth and discontinuous features of an exact solution to the linear advection equation 
	\begin{equation*}
		\frac{\partial u}{\partial t}+\nabla\cdot(\mathbf{v}u)=0\quad\mbox{in}\ \Omega=(0,1)^2
	\end{equation*}
	with the velocity field $\mathbf{v}(x,y)=(0.5-y,x-0.5)^\top$. The initial and boundary conditions are taken from \cite{kuzmin2020,leveque1996}. Since the exact solution is periodic in time, it coincides with the initial condition after each full revolution $t\in \mathbb{N}$, and we stop computations at the final time $t=1$.
	
	We use the mesh size $h$ corresponding to $N_h=257^2$ DoFs for a given polynomial degree $p$, and employ continuous finite elements. The results are shown in Fig.~\ref{fig:sbr}. For linear and quadratic polynomials, the new CV-WENO variant produces solutions comparable to classical CC-WENO schemes. For quartic Lagrange polynomials, however, the CV-WENO scheme reproduces both smooth and discontinuous solution features far more accurately, while the CC-WENO scheme exhibits a clear loss of accuracy. This behavior has previously been reported in \cite{kuzmin2023a} and the reason for this is twofold. First, the CC-WENO reconstruction only uses information from von Neumann neighbors and excludes vertex neighbors, resulting in a significant loss of spatial information that becomes increasingly severe for higher polynomial degrees and coarser meshes. Second, the CC-WENO scheme extrapolates high-order polynomials from neighboring elements into the target element, a procedure that becomes increasingly ill-posed as the polynomial degree increases. The CV-WENO scheme resolves both issues.
	
	As shown in Table~\ref{tab:sbr}, equipping the CV-WENO scheme with shock-capturing WENO quadrature reduces the $L^1$ errors compared to the standard CV-WENO formulation.
	
	\begin{figure}[h!]
		\centering
		\small
		
		% First row
		\begin{minipage}[t]{0.32\textwidth}
			\centering
			(a) CC-WENO, $p=1$,\\$u_h\in[-0.06179,1.08641]$ \\[0.2cm] 
			\includegraphics[width=0.95\textwidth,trim=0 0 0 0,clip]{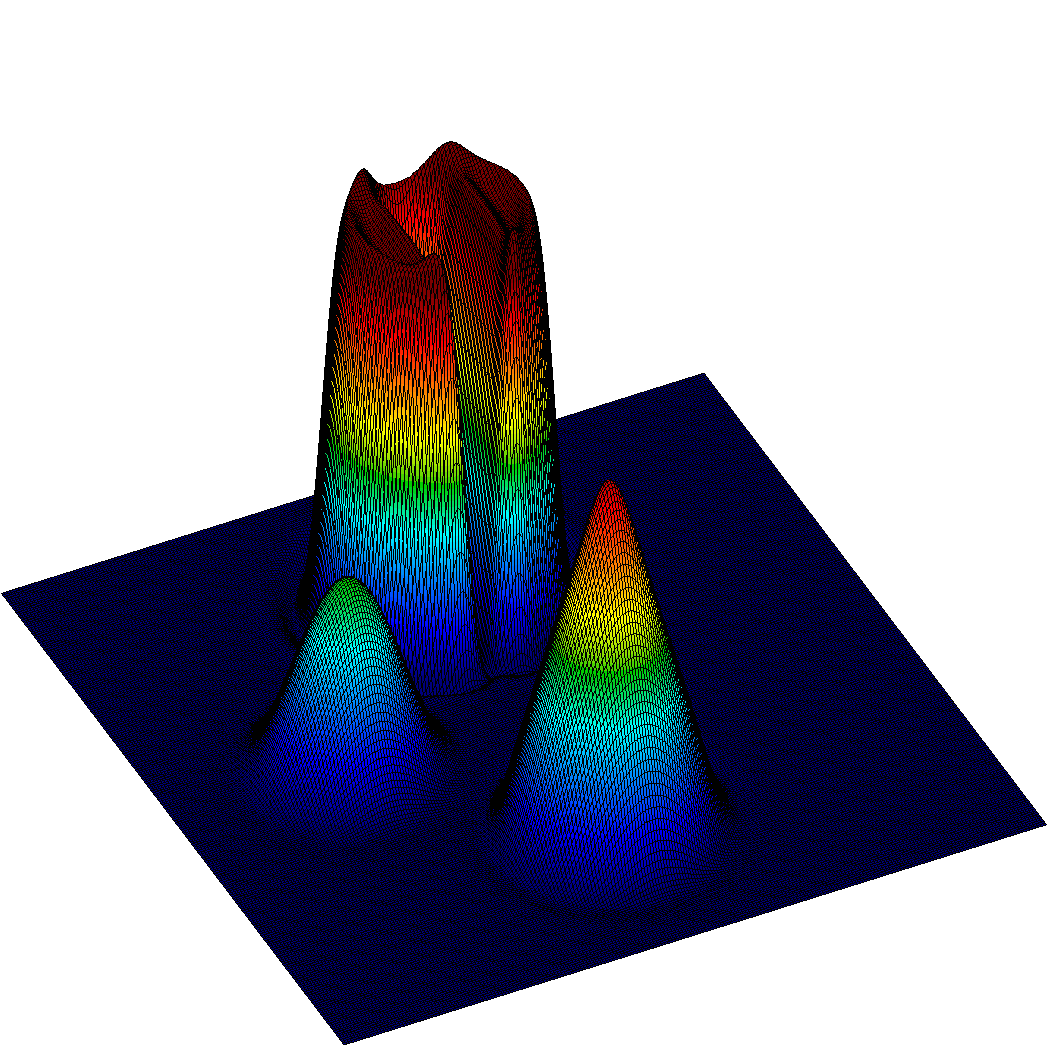}
		\end{minipage}%
		\hfill
		\begin{minipage}[t]{0.32\textwidth}
			\centering
			(b) CV-WENO, $p=1$,\\$u_h\in[-0.06710,1.08863]$\\[0.2cm]
			\includegraphics[width=0.95\textwidth,trim=0 0 0 0,clip]{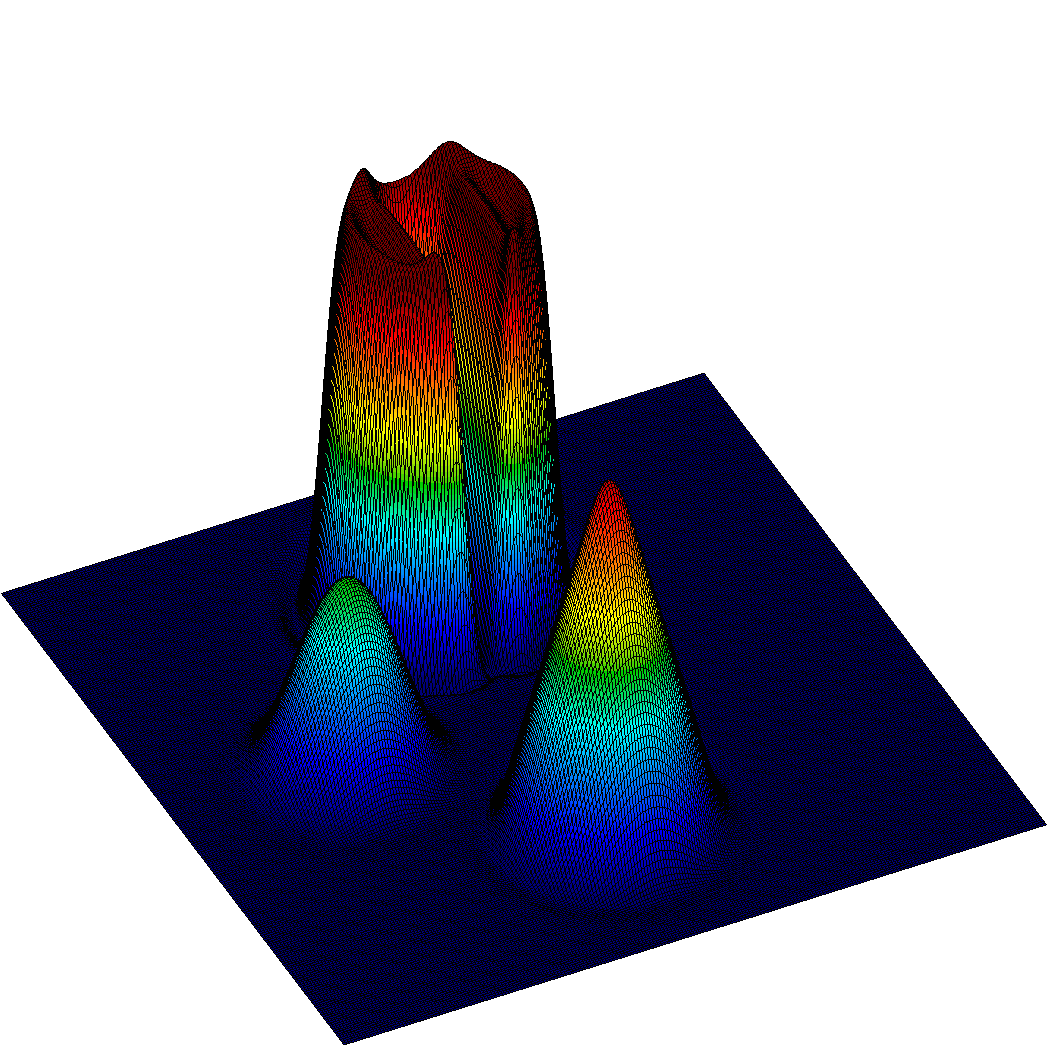}
		\end{minipage}%
		\hfill
		\begin{minipage}[t]{0.32\textwidth}
			\centering
			(c) CV-WENO-SC, $p=1$,\\$u_h\in[-0.06711,1.08865]$\\[0.2cm]
			\includegraphics[width=0.95\textwidth,trim=0 0 0 0,clip]{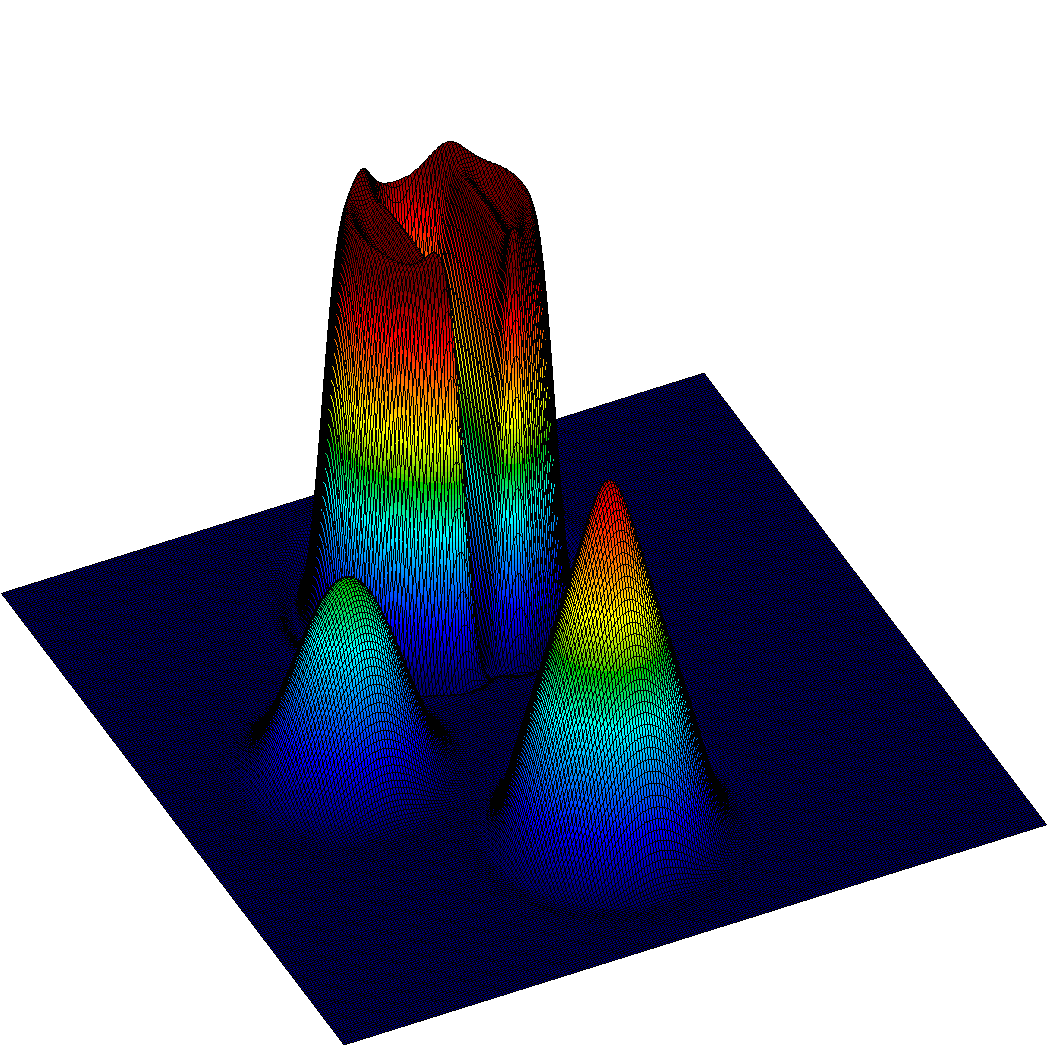}
		\end{minipage}
		
		\vskip0.5cm % vertical space between rows
		
		% Second row
		\begin{minipage}[t]{0.32\textwidth}
			\centering
			(d) CC-WENO, $p=2$,\\$u_h\in[-0.00439,1.00751]$ \\[0.2cm] 
			\includegraphics[width=0.95\textwidth,trim=0 0 0 0,clip]{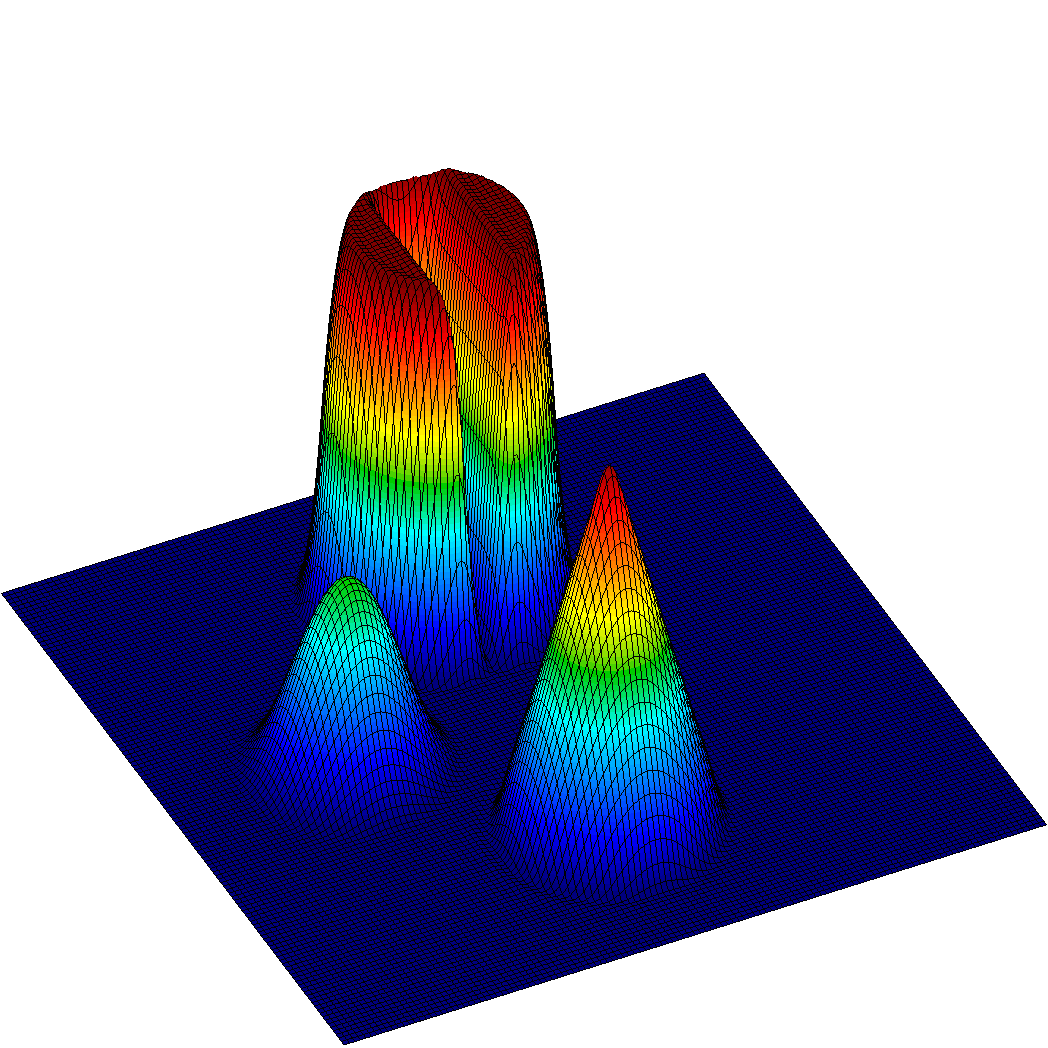}
		\end{minipage}%
		\hfill
		\begin{minipage}[t]{0.32\textwidth}
			\centering
			(e) CV-WENO, $p=2$,\\$u_h\in[-0.00975,1.00819]$\\[0.2cm]
			\includegraphics[width=0.95\textwidth,trim=0 0 0 0,clip]{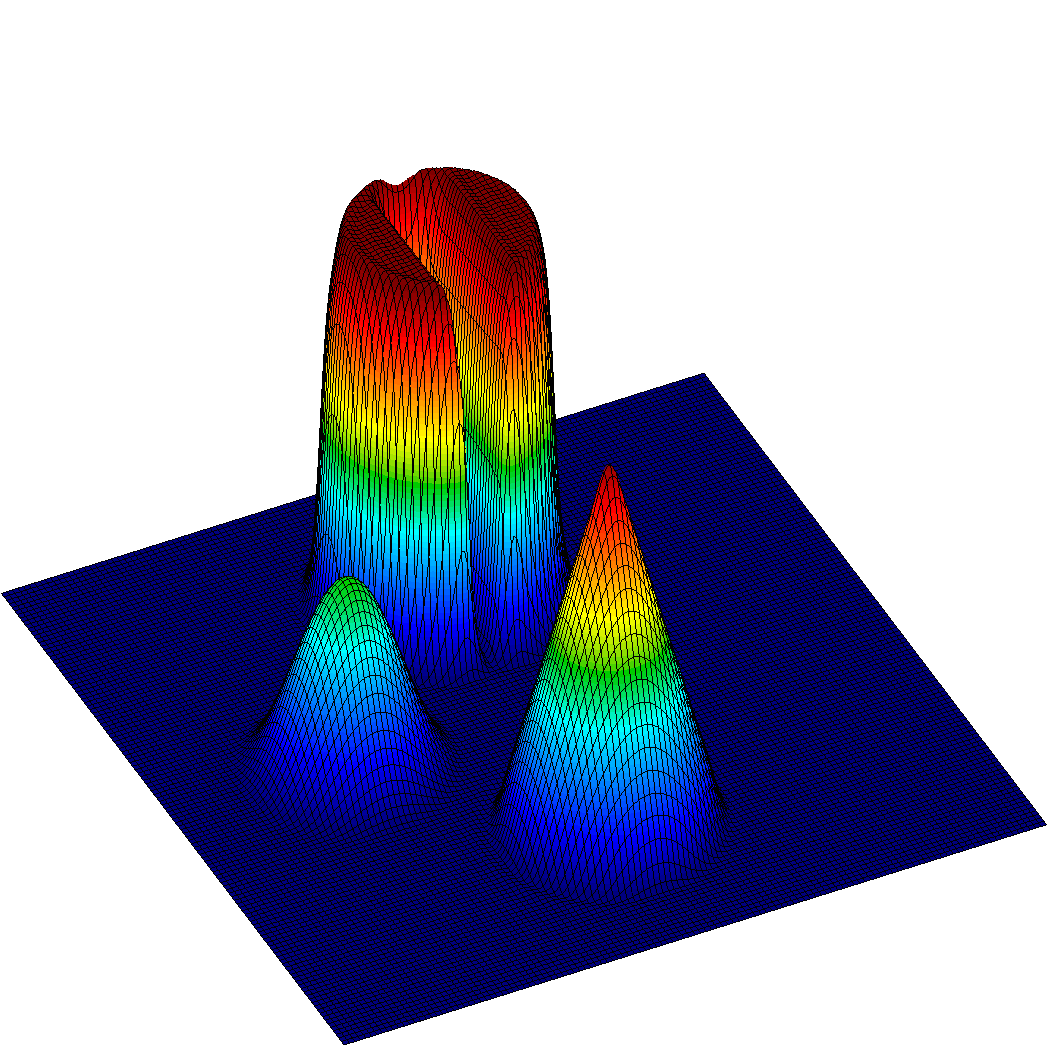}
		\end{minipage}%
		\hfill
		\begin{minipage}[t]{0.32\textwidth}
			\centering
			(f) CV-WENO-SC, $p=2$,\\$u_h\in[-0.00984,1.00832]$\\[0.2cm]
			\includegraphics[width=0.95\textwidth,trim=0 0 0 0,clip]{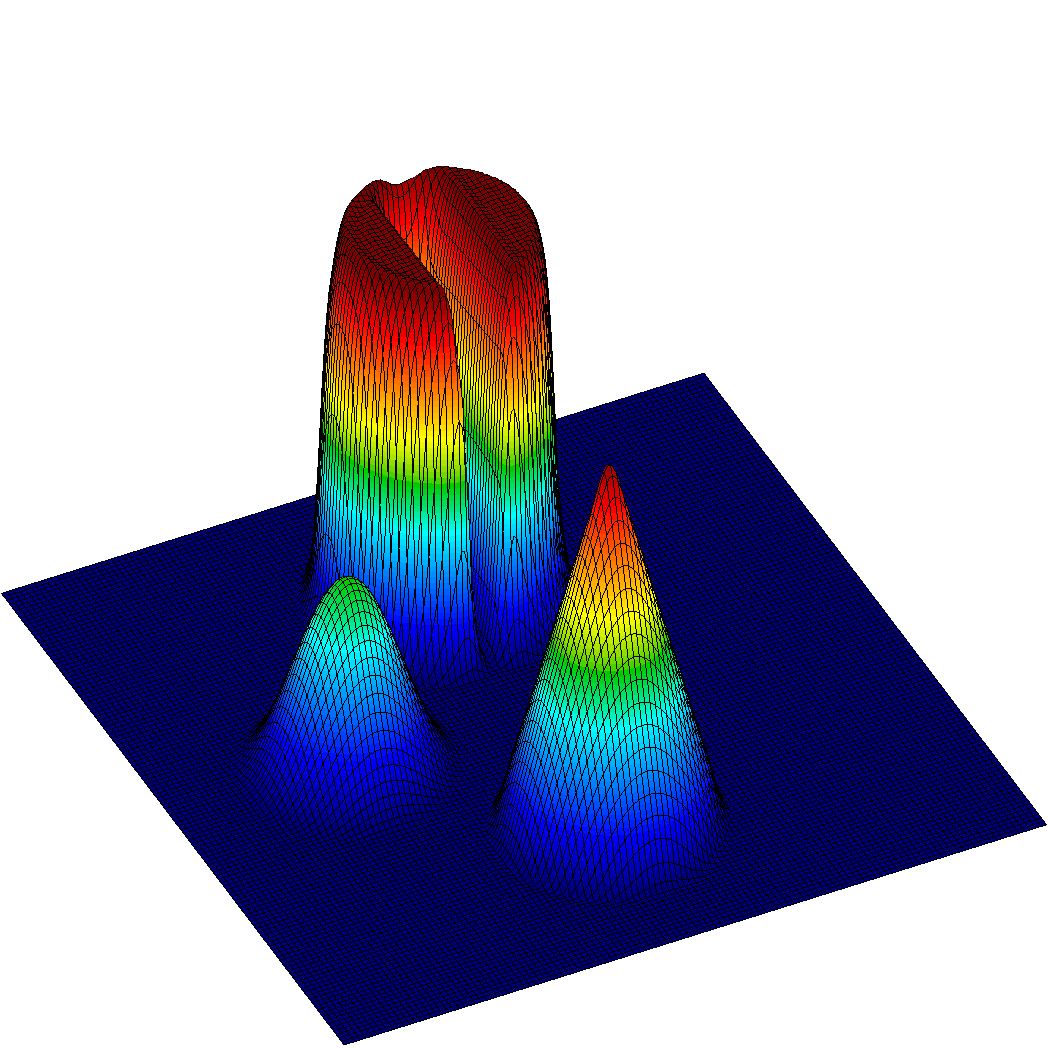}
		\end{minipage}
		
		\vskip0.5cm % vertical space between rows
		
		\begin{minipage}[t]{0.32\textwidth}
			\centering
			(g) CC-WENO, $p=4$,\\$u_h\in[-1.7e-05,0.85188]$ \\[0.2cm] 
			\includegraphics[width=0.95\textwidth,trim=0 0 0 0,clip]{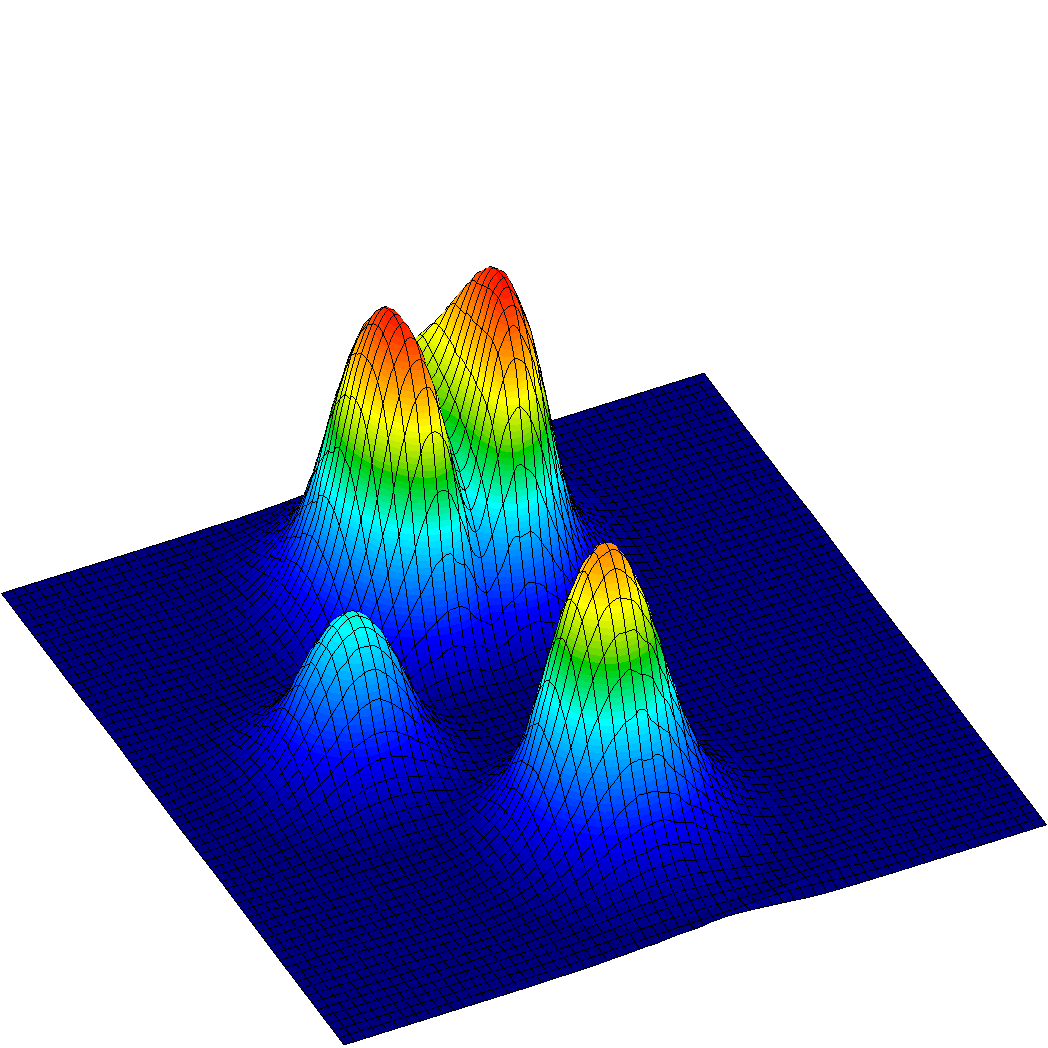}
		\end{minipage}%
		\hfill
		\begin{minipage}[t]{0.32\textwidth}
			\centering
			(h) CV-WENO, $p=4$,\\$u_h\in[-0.00310,1.00466]$\\[0.2cm]
			\includegraphics[width=0.95\textwidth,trim=0 0 0 0,clip]{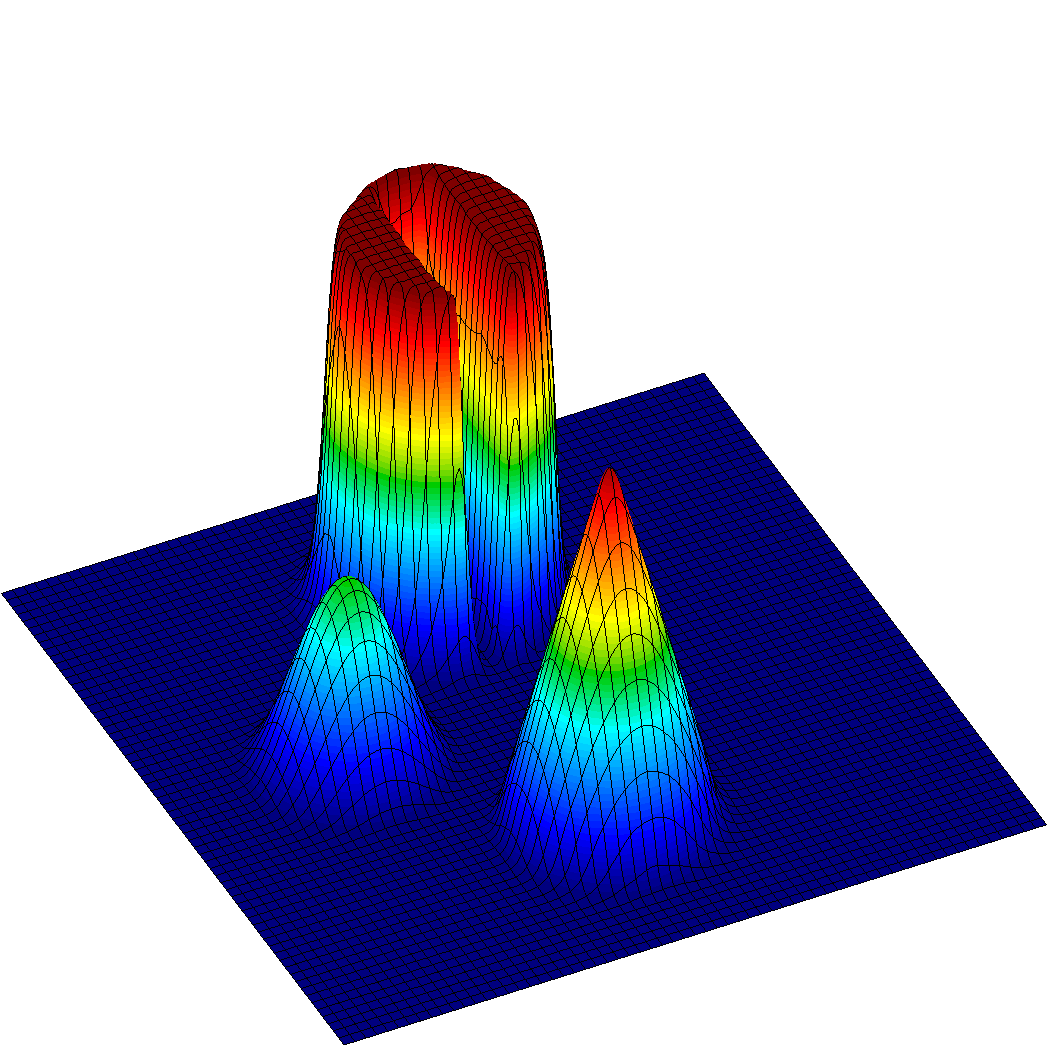}
		\end{minipage}%
		\hfill
		\begin{minipage}[t]{0.32\textwidth}
			\centering
			(i) CV-WENO-SC, $p=4$,\\$u_h\in[-0.05656,1.02060]$\\[0.2cm]
			\includegraphics[width=0.95\textwidth,trim=0 0 0 0,clip]{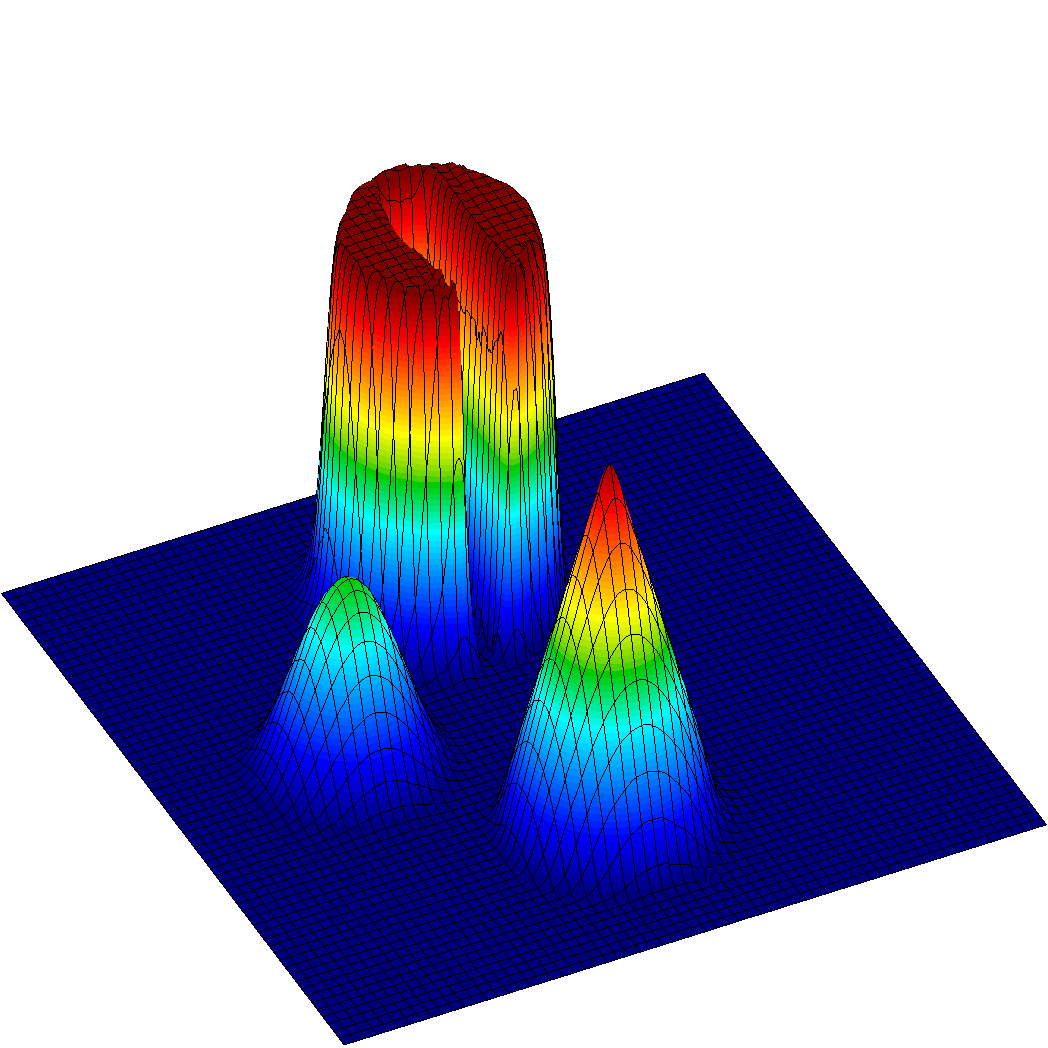}
		\end{minipage}
		
		\caption{Numerical solutions to the solid body rotation problem \cite{leveque1996} at $t=1$ obtained using $N_h=257^2$ and $p\in\{1,2,4\}$.}
		\label{fig:sbr}
	\end{figure}
	
	\begin{table}[htb!]
		\centering
		\begin{tabular}{cccc}
			\hline
			 & $p=1$ & $p=2$ & $p=4$ \\
			\hline
			CC-WENO    & 1.60595e-02 & 1.57469e-02 & 5.43568e-02 \\
			CV-WENO    & 1.59866e-02 & 1.18561e-02 & 1.17923e-02 \\
			CV-WENO-SC & 1.59861e-02 & 1.17689e-02 & 1.08783e-02 \\
			\hline
		\end{tabular}
		\caption{Errors $\|u_h-u_{\textrm{ex}}\|_{L^1(\Omega)}$ to the solid body rotation problem \cite{leveque1996} at $t=1$ obtained using $N_h=257^2$ and $p\in\{1,2,4\}$.}
		\label{tab:sbr}
	\end{table}
\end{subsection}

\begin{subsection}{One-dimensional inviscid Burgers equation}
	In the next experiment, we solve the one-dimensional inviscid Burgers equation
	\begin{equation*}
		\frac{\partial u}{\partial t}+\frac{\partial (u^2/2)}{\partial x}=0\quad\mbox{in}\ \Omega=(0,1).
	\end{equation*}
	The initial condition
	\begin{equation}
		u_0(x)=\sin(2\pi x)
		\label{eq:burinitial}
	\end{equation}
	remains smooth up to the critical time $t_c=\frac{1}{2\pi}$ of shock formation. We perform grid convergence studies at $t=0.1$, for which the solution is sufficiently smooth, and report the results in Tables~\ref{tab:bur1} and \ref{tab:bur2}. The time step is chosen small enough to make the spatial error dominant. As before, we use continuous Lagrange finite elements. Since the shock-capturing WENO quadrature is inactive in smooth regions, the numerical errors produced by the CV-WENO and CV-WENO-SC scheme are identical and therefore reported only once. It is observed that all methods achieve optimal convergence rates. Moreover, consistent with the behavior observed for the linear solid body rotation benchmark, CV-WENO outperforms CC-WENO on coarse meshes and reaches the asymptotic regime faster.
	
	\begin{table}[t!]
		\small
		\centering
		\begin{tabular}{ccccccc}
			\hline
			$E_h$ & $p=1$ & EOC & $p=2$ & EOC & $p=3$ & EOC \\
			\hline
			48	 & 1.76E-03 &      & 5.30E-05 &      & 1.43E-05 &	   \\
			64   & 7.73E-04 & 2.86 & 2.35E-05 & 2.83 & 7.27E-06 & 2.35 \\
			96   & 2.72E-04 & 2.58 & 8.21E-06 & 2.59 & 2.93E-06 & 2.24 \\
			128	 & 1.45E-04 & 2.19 & 3.87E-06 & 2.61 & 1.50E-06 & 2.33 \\
			192	 & 6.50E-05 & 1.98 & 1.21E-06 & 2.87 & 4.79E-07 & 2.82 \\
			256	 & 3.67E-05 & 1.99 & 4.82E-07 & 3.20 & 1.76E-07 & 3.48 \\
			384  & 1.64E-05 & 1.99 & 1.21E-07 & 3.41 & 3.75E-08 & 3.81 \\
			512	 & 9.23E-06 & 2.00 & 4.59E-08 & 3.37 & 1.21E-08 & 3.93 \\
			768	 & 4.11E-06 & 2.00 & 1.22E-08 & 3.27 & 2.46E-09 & 3.93 \\
			1024 & 2.31E-06 & 2.00 & 4.87E-09 & 3.19 &          &	   \\
			1536 & 1.03E-06 & 1.99 &          &      &          &	   \\
			\hline	
		\end{tabular}
		\caption{$L^1(\Omega)$ error norms and convergence rates for the CC-WENO scheme applied to the 1D inviscid Burgers equation with initial condition \eqref{eq:burinitial} at $t=0.1$.}
		\label{tab:bur1}
	\end{table}
	\begin{table}[htp!]
		\small
		\centering
		\begin{tabular}{ccccccc}
			\hline
			$E_h$ & $p=1$ & EOC & $p=2$ & EOC & $p=3$ & EOC \\
			\hline
			48	 & 1.83E-03 &      & 4.70E-05 &      & 4.03E-06 &	   \\
			64   & 7.59E-04 & 3.06 & 1.91E-05 & 3.13 & 1.33E-06 & 3.85 \\
			96   & 2.66E-04 & 2.59 & 5.54E-06 & 3.05 & 2.73E-07 & 3.91 \\
			128	 & 1.45E-04 & 2.11 & 2.32E-06 & 3.03 & 8.84E-08 & 3.92 \\
			192	 & 6.50E-05 & 1.98 & 6.81E-07 & 3.02 & 1.80E-08 & 3.93 \\
			256	 & 3.67E-05 & 1.99 & 2.87E-07 & 3.00 & 5.85E-09 & 3.91 \\
			384  & 1.64E-05 & 1.99 & 8.48E-08 & 3.01 & 1.21E-09 & 3.89 \\
			512	 & 9.22E-06 & 2.00 & 3.57E-08 & 3.01 & 3.96E-10 & 3.88 \\
			768	 & 4.11E-06 & 1.99 & 1.06E-08 & 2.99 & 8.29E-11 & 3.86 \\
			1024 & 2.31E-06 & 2.00 & 4.46E-09 & 3.01 &          &	   \\
			1536 & 1.03E-06 & 1.99 &          &      &          &	   \\
			\hline
		\end{tabular}
		\caption{$L^1(\Omega)$ error norms and convergence rates for the CV-WENO(-SC) scheme applied to the 1D inviscid Burgers equation with initial condition \eqref{eq:burinitial} at $t=0.1$.}
		\label{tab:bur2}
	\end{table}
	
	\begin{remark}
		Similar convergence behavior is observed for the one-dimensional linear advection equation with a smoothed step function, taken from \cite[Sec.~7.1]{hajduk2020b}, using a constant velocity $v=1$ (not shown here). As before, all schemes converge at the optimal rate $p+1$. The DG variant of the proposed methods also exhibits optimal convergence at the same rate.
	\end{remark}
	
	By extending the final time to $t=1>t_c$, after the shock has fully developed, we assess the shock-capturing properties of the schemes. With the total number of degrees of freedom fixed at $N_h=128$, the results in Fig.~\ref{fig:bur} show that all WENO schemes produce nonoscillatory solutions with nearly identical profiles and increasing the polynomial degree $p$ has little effect. Due to the self-steepening nature of shock fronts, relatively diffusive schemes perform well in this benchmark. Even the uniform choice $\gamma_e=1$ in \eqref{eq:stab}, which corresponds to a low-order Lax-Friedrichs-type scheme, yields satisfactory results, as reported in \cite{kuzmin2023a,vedral-arxiv}.

	\begin{figure}[t!]
		\centering
		\small

		\begin{minipage}[t]{0.32\textwidth}
			\centering
			(a) CC-WENO\\[0.2cm] 
			\includegraphics[width=0.95\textwidth,trim=0 0 0 0,clip]{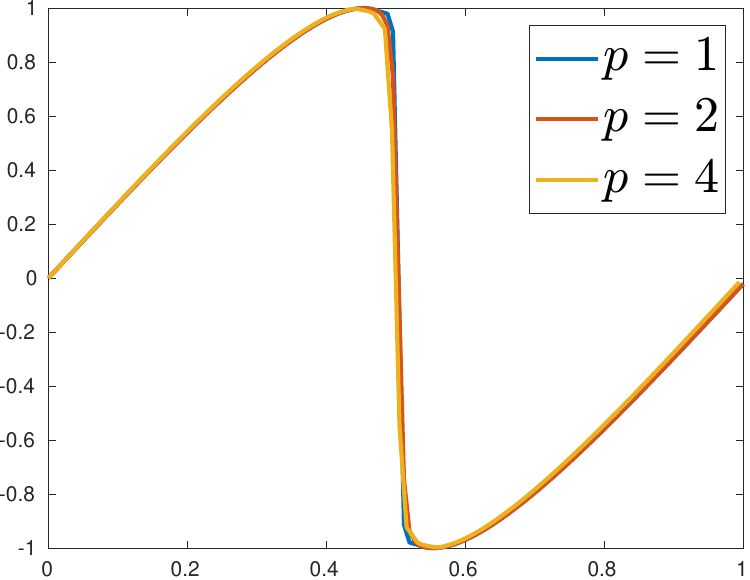}
		\end{minipage}%
		\hfill
		\begin{minipage}[t]{0.32\textwidth}
			\centering
			(b) CV-WENO\\[0.2cm]
			\includegraphics[width=0.95\textwidth,trim=0 0 0 0,clip]{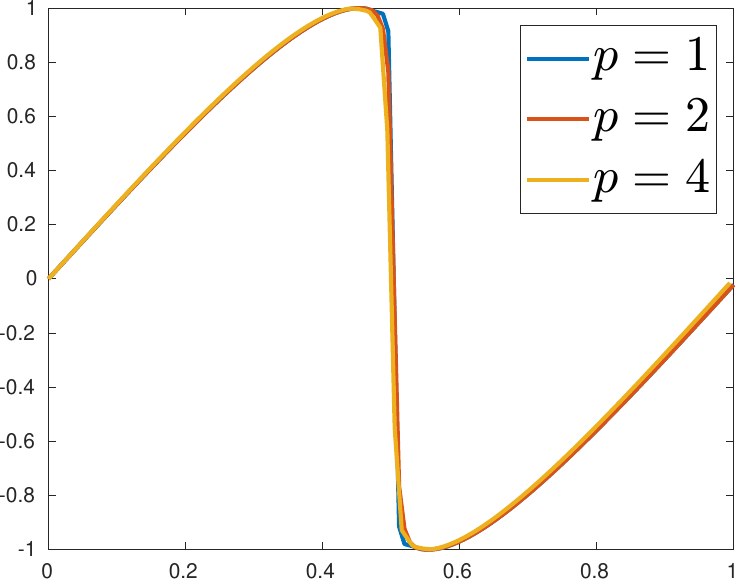}
		\end{minipage}%
		\hfill
		\begin{minipage}[t]{0.32\textwidth}
			\centering
			(c) CV-WENO-SC\\[0.2cm]
			\includegraphics[width=0.95\textwidth,trim=0 0 0 0,clip]{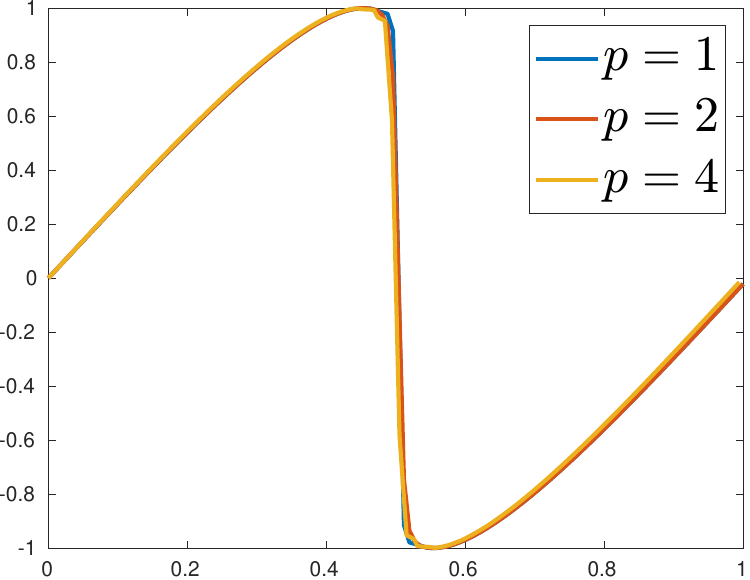}
		\end{minipage}
		
		\caption{Numerical solutions to the 1D inviscid Burgers equation at $t=0.2$ obtained using $N_h=128$ and $p\in\{1,2,4\}$.}
		\label{fig:bur}
	\end{figure}
\end{subsection}

\begin{subsection}{KPP problem}
	In the next experiment, we consider the KPP problem \cite{kurganov2007}, a challenging nonlinear test problem for assessing entropy stability of high-order numerical schemes. We solve \eqref{eq:pde} with the nonconvex flux function $\mathbf{f}(u)=(\sin(u),\cos(u))$ on the computational domain $\Omega=(-2,2)\times(-2.5,1.5)$. The initial condition
	\begin{align*}
		u_0(x,y)=\begin{cases}
			\frac{7\pi}{2} & \mbox{if}\quad \sqrt{x^2+y^2}\le 1,\\
			\frac{\pi}{4} & \mbox{otherwise},
		\end{cases}
	\end{align*}
	generates an exact entropy solution exhibiting a rotational wave structure. Standard schemes often fail to converge to the exact entropy solution without an appropriate entropy fix. For dissipation-based methods, introducing sufficient dissipation near shocks is essential to converge to the correct entropy solution. A global upper bound for the maximum wave speed, required to compute the viscosity parameter $\nu_e$ in \eqref{eq:lostab}, is given by $\lambda_e=1.0$. More accurate bounds can be found in \cite{guermond2017}.
	
	Simulations are performed on a uniform mesh with $N_h=257^2$ DoFs using continuous Lagrange finite elements of order $p=2$. The results, shown in Fig.~\ref{fig:kpp}, indicate that the CV-WENO scheme produces nonoscillatory solutions that accurately capture the spiral wave structure, with accuracy comparable to the CC-WENO variant. The results for $p=1$ and $p=4$ (not shown here) are similar. Although a (semi-)discrete entropy inequality cannot be formally proven for our scheme, as is the case for most WENO schemes reported in the literature, the results suggest that it is at least numerically entropy stable.
	\begin{figure}[t!]
		\centering
		\small
		
		% First row
		\begin{minipage}[t]{0.32\textwidth}
			\centering
			(a) CC-WENO,\\$u_h\in[0.77223,11.00050]$ \\[0.2cm] 
			\includegraphics[width=0.95\textwidth,trim=0 0 0 0,clip]{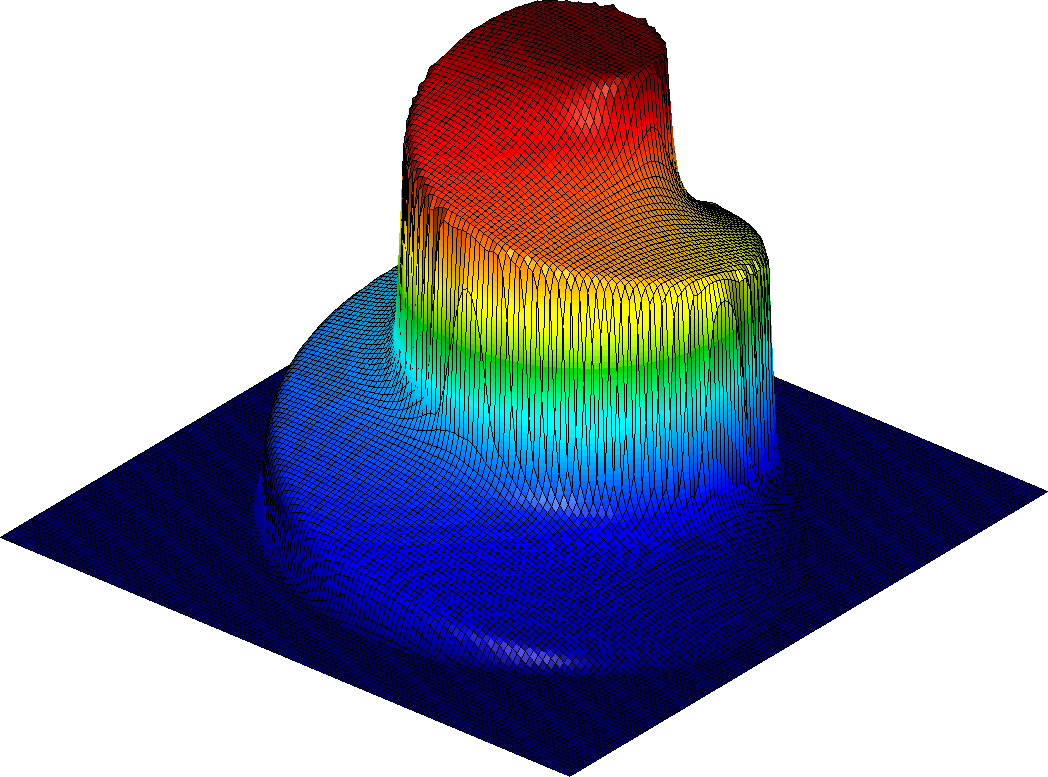}
		\end{minipage}%
		\hfill
		\begin{minipage}[t]{0.32\textwidth}
			\centering
			(b) CV-WENO,\\$u_h\in[0.77155,10.99950]$\\[0.2cm]
			\includegraphics[width=0.95\textwidth,trim=0 0 0 0,clip]{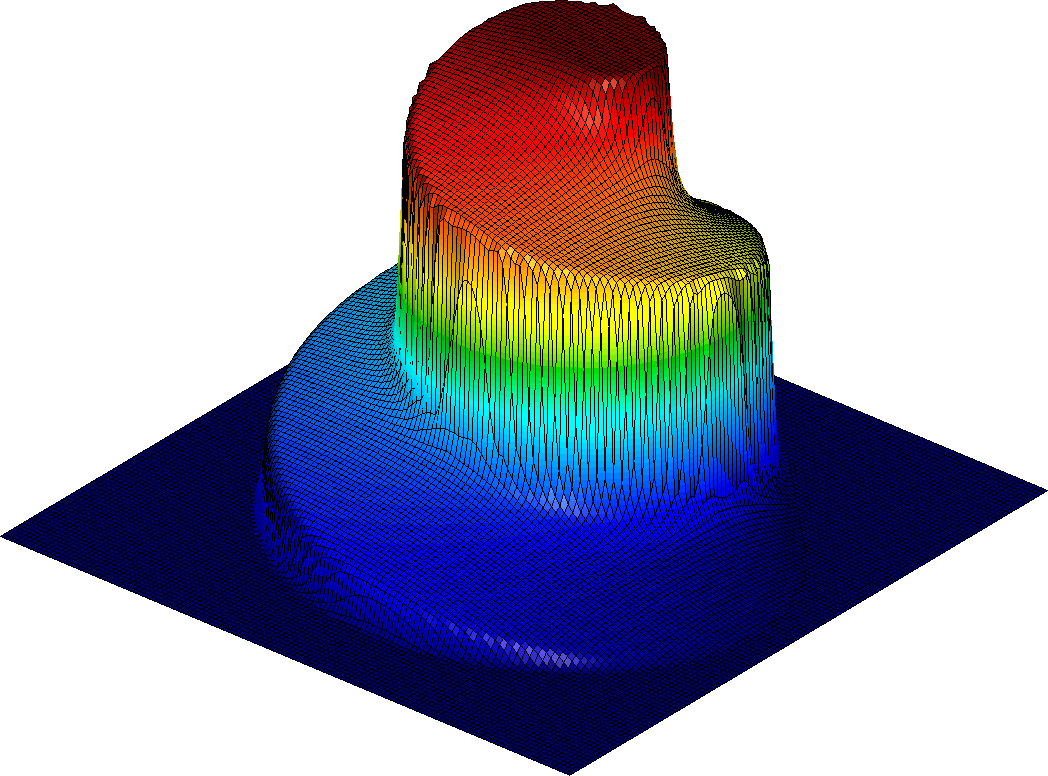}
		\end{minipage}%
		\hfill
		\begin{minipage}[t]{0.32\textwidth}
			\centering
			(c) CV-WENO-SC,\\$u_h\in[0.68227,11.02270]$\\[0.2cm]
			\includegraphics[width=0.95\textwidth,trim=0 0 0 0,clip]{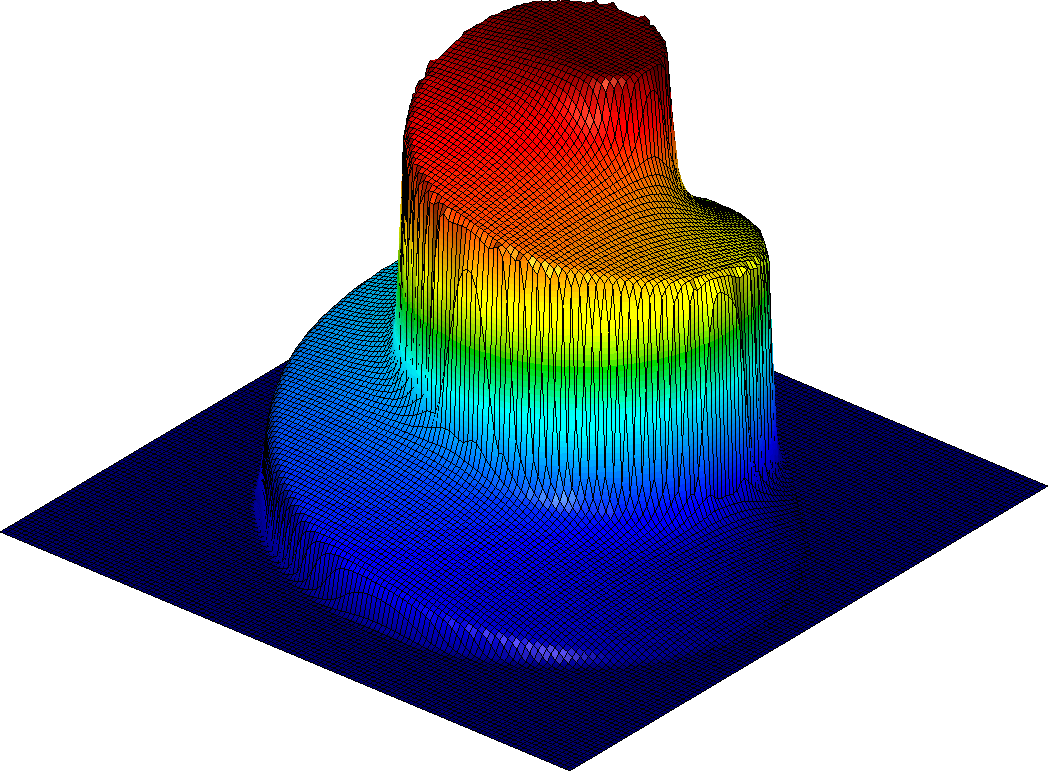}
		\end{minipage}
		
		\vskip0.25cm % vertical space between rows
		
		% Second row
		\begin{minipage}[t]{0.32\textwidth}
			\centering
			%(d) CC-WENO, top view,\\$u_h\in[-0.00439,1.00751]$ \\[0.2cm] 
			\includegraphics[width=0.95\textwidth,trim=0 0 0 0,clip]{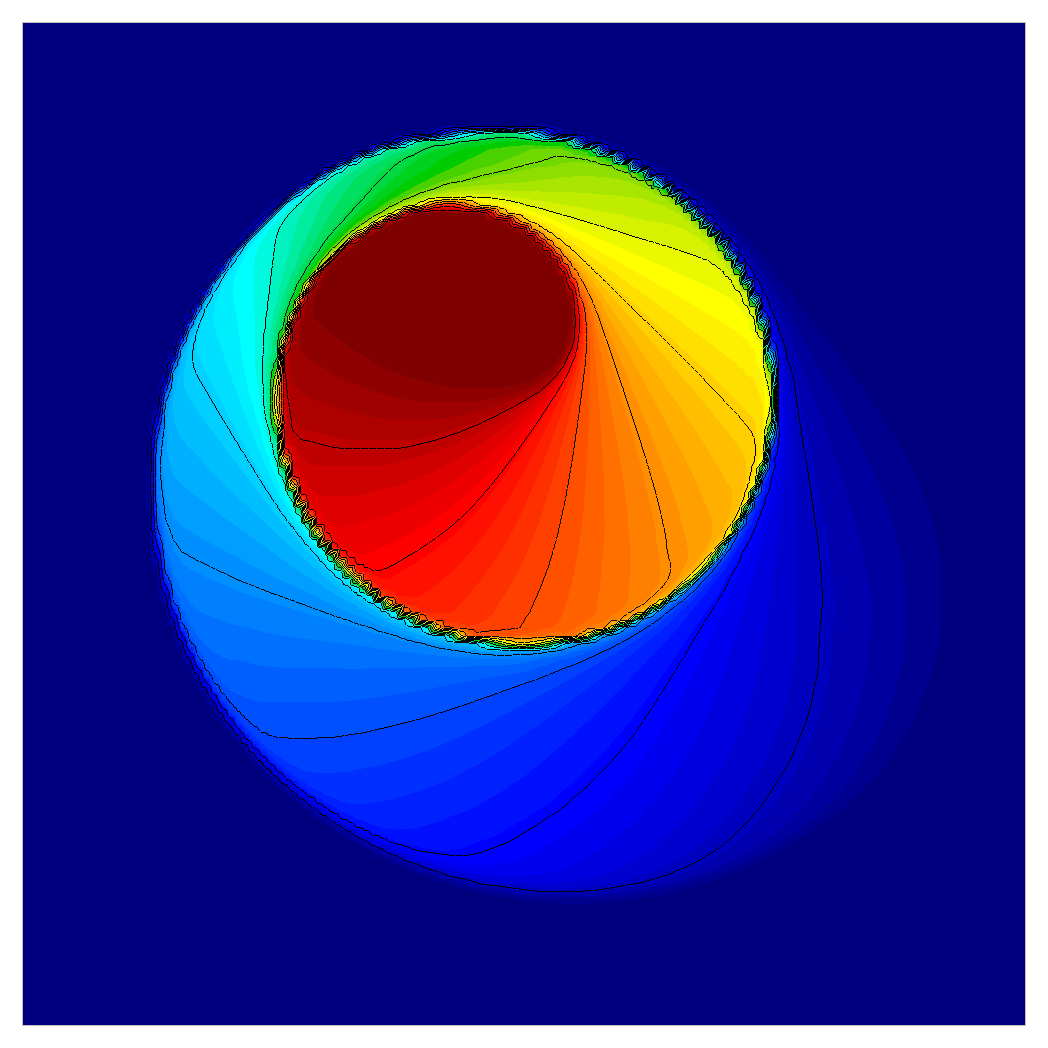}
		\end{minipage}%
		\hfill
		\begin{minipage}[t]{0.32\textwidth}
			\centering
			%(e) CV-WENO, top view,\\$u_h\in[-0.00975,1.00819]$\\[0.2cm]
			\includegraphics[width=0.95\textwidth,trim=0 0 0 0,clip]{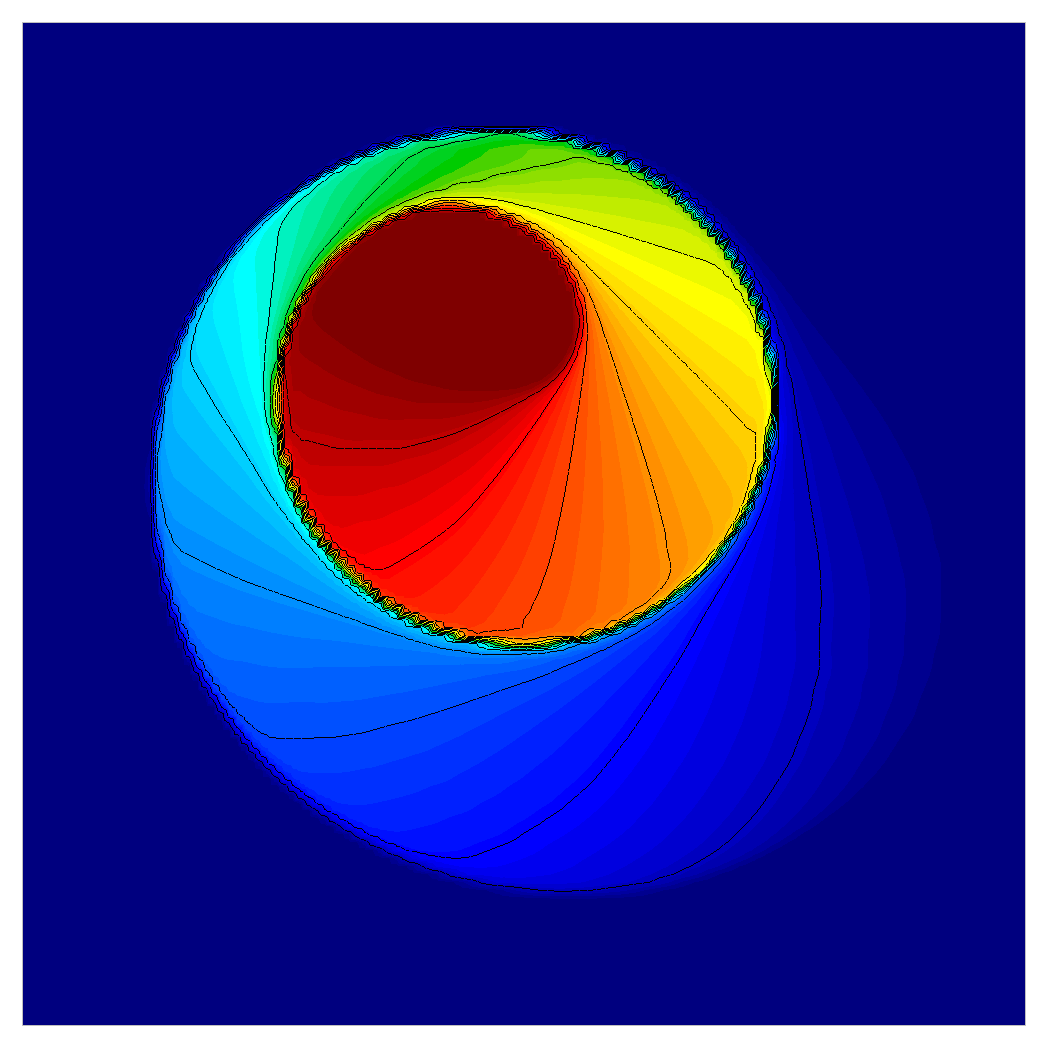}
		\end{minipage}%
		\hfill
		\begin{minipage}[t]{0.32\textwidth}
			\centering
			%(f) CV-WENO-SC, top view,\\$u_h\in[-0.00984,1.00832]$\\[0.2cm]
			\includegraphics[width=0.95\textwidth,trim=0 0 0 0,clip]{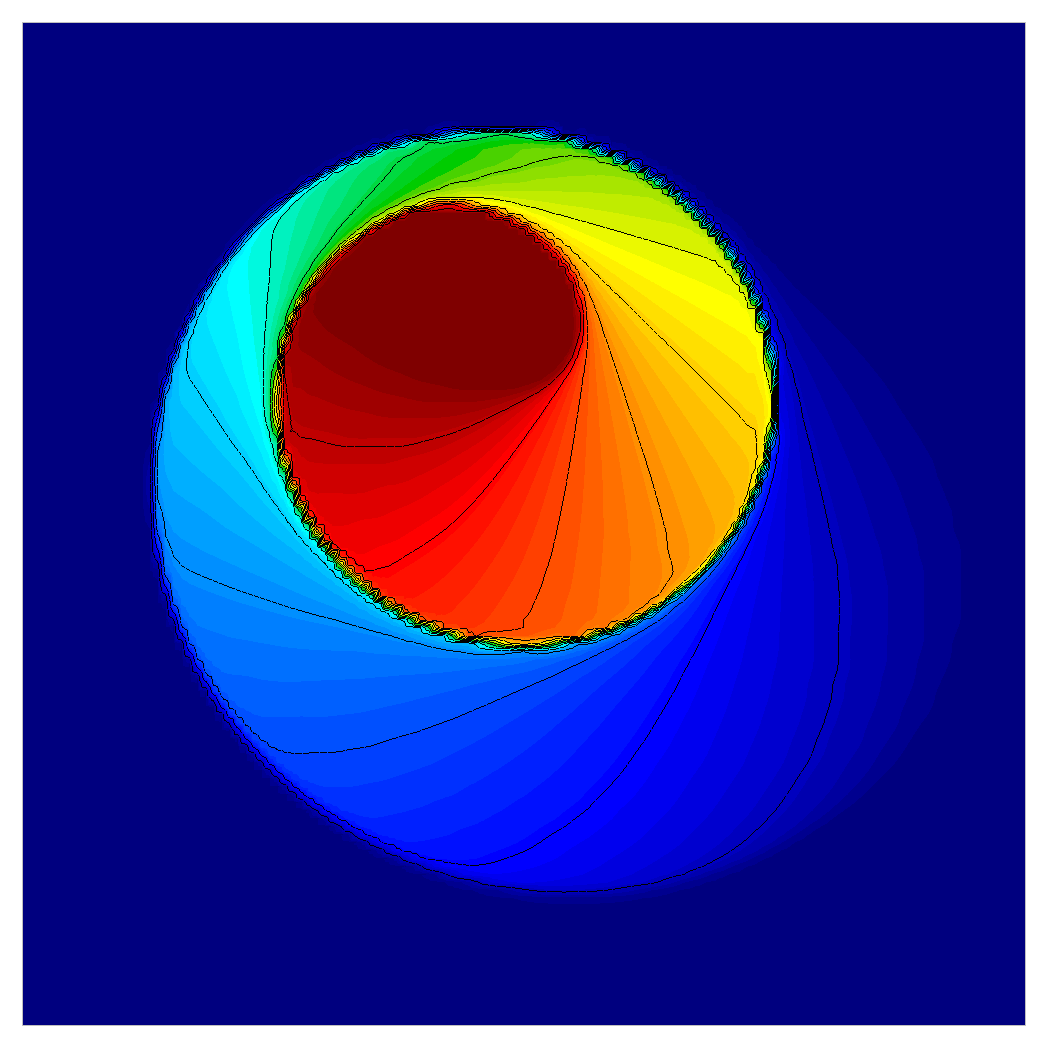}
		\end{minipage}
		
		\caption{Numerical solutions to the KPP problem \cite{kurganov2007} at $t=1$ obtained using $N_h=257^2$ and $p=2$.}
		\label{fig:kpp}
	\end{figure}
\end{subsection}
\begin{subsection}{Euler equations of gas dynamics}
	We consider the Euler equations of gas dynamics which describe the conservation of mass, momentum and total energy. The vector of conserved variables and the corresponding flux tensor in \eqref{eq:pde} are given by
	\begin{equation*}
		U = \begin{bmatrix}
			\varrho \\
			\varrho \mathbf{v}\\
			\varrho E
		\end{bmatrix}\in \mathbb{R}^{d+2}, \quad 
		\mathbf{F(U)}=\begin{bmatrix}
			\varrho \mathbf{v} \\
			\varrho \mathbf{v} \bigotimes \mathbf{v}+p\mathbf{I}\\
			(\varrho E + p)\mathbf{v}
		\end{bmatrix}\in \mathbb{R}^{(d+2) \times d},
	\end{equation*}
	where $\varrho$, $\mathbf{v}$, $E$ denote the density, velocity and specific total energy, respectively. The identity matrix is denoted by $\mathbf{I}$. The pressure $p$ is calculated using the polytropic ideal gas law
	\begin{equation*}
		p=\varrho e (\gamma-1),
	\end{equation*}
	where $\varrho e$ and $\gamma$ denote the internal energy and the heat capacity ratio, respectively. In all numerical experiments, we set $\gamma=1.4$.
	
	Following the approach in \cite{kuzmin2023a,vedral-arxiv}, the blending factor $\gamma_e$ in \eqref{eq:gamma} is computed from the density field and subsequently applied to all variables, yielding significant computational savings without compromising robustness, as reported in \cite{pirozzoli2002}. 
	
	In what follows, we do not present results obtained with the CC-WENO scheme. Numerical solutions of the Euler equations using CC-WENO are available in \cite{vedral-arxiv} for the DG formulation and in the recent work \cite{kuzmin2025a} for the CG formulation.
\end{subsection}
\begin{subsubsection}{Modified Sod shock tube}
	We assess the entropy stability of our scheme for hyperbolic systems using the modified Sod shock tube problem \cite{toro2013}, a modification of the classical Sod shock tube \cite{sod1978}. The initial conditions are prescribed as
	\begin{align*}
		(\rho_0,v_0,p_0)(x)=\begin{cases}
			(1,0.75,1) &\text{if } x < 0.25,\\
			(0.125,0,0.1)&\text{otherwise},
		\end{cases}
	\end{align*}
	within the computational domain $\Omega=(0,1)$. The left boundary at $x=0$ acts as a supersonic inlet, while the right boundary at $x=1$ is reflective. The initial discontinuity at $x=0.25$ generates a rarefaction wave, a contact discontinuity, and a shock, all of which the numerical scheme should resolve accurately. High-order methods are known to produce entropy glitches near the sonic point of the rarefaction wave, highlighting the importance of robust entropy-stable methods.
	
	We run simulations up to the final time $t=0.2$ on uniform meshes. Although DG finite elements are commonly regarded as better suited for hyperbolic problems than CG finite elements, we investigate this assumption through our numerical tests. Both approaches are compared using roughly the same number of degrees of freedom per variable ($N_{h,1}=257$ for CG and $N_{h,2}=256$ for DG) to ensure a fair evaluation. We further aim to demonstrate that our scheme can handle very high-order polynomials.
	
	The numerical results obtained with the CV-WENO and CV-WENO-SC variants, using continuous and discontinuous finite elements, are shown in Figs.~\ref{fig:msod1} and \ref{fig:msod2}, respectively. The results indicate that the CG formulation of the CV-WENO strategy performs comparably to its DG counterpart in terms of accuracy and shock-capturing properties. Remarkably, we obtain reasonable solutions even for polynomial degrees as high as $p=32$ on only eight elements. The shock-capturing quadrature produces less smearing near the rarefaction wave. Moreover, for all polynomial degrees and for both CG and DG formulations, the numerical solutions remain entropy stable, eliminating the need for additional (semi-)discrete entropy fixes. These findings are consistent with the observations from the KPP problem.
	
	\begin{figure}[t!]
		\centering
		\small
		\begin{minipage}[t]{0.48\textwidth}
			\centering
			(a) CV-WENO\\[0.2cm] 
			\includegraphics[width=0.95\textwidth,trim=0 0 0 0,clip]{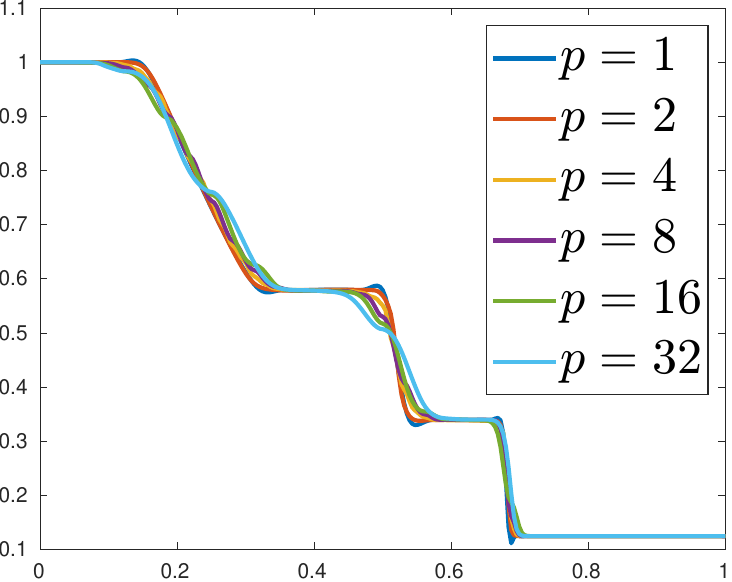}
		\end{minipage}%
		\hfill
		\begin{minipage}[t]{0.48\textwidth}
			\centering
			(b) CV-WENO-SC\\[0.2cm]
			\includegraphics[width=0.95\textwidth,trim=0 0 0 0,clip]{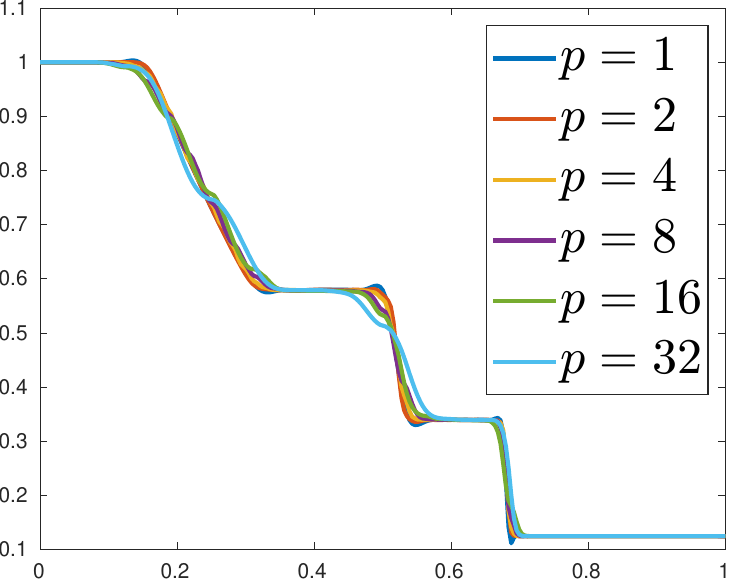}
		\end{minipage}%
		\caption{Modified Sod shock tube, density profiles $\varrho$ at $t=0.2$ obtained using continuous finite elements with $N_h=257$ DoFs per variable.}
		\label{fig:msod1}
	\end{figure} 
	\begin{figure}[htp!]
		\centering
		\small
		\begin{minipage}[t]{0.48\textwidth}
			\centering
			(a) CV-WENO\\[0.2cm] 
			\includegraphics[width=0.95\textwidth,trim=0 0 0 0,clip]{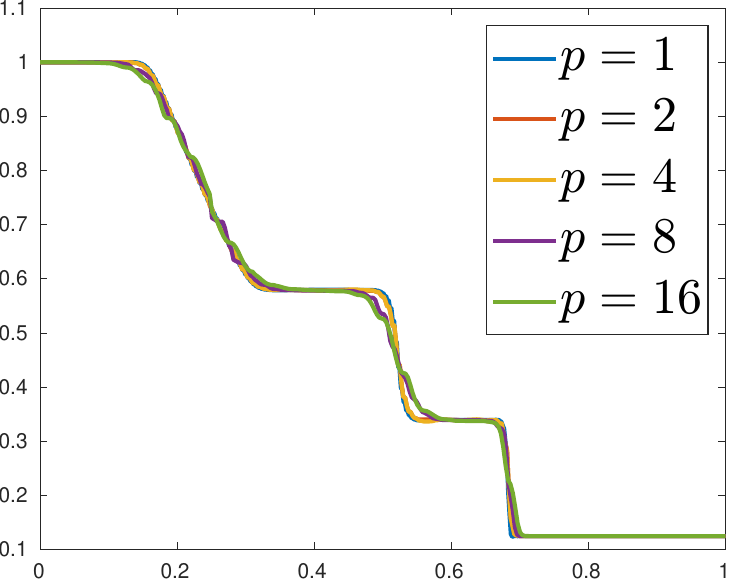}
		\end{minipage}%
		\hfill
		\begin{minipage}[t]{0.48\textwidth}
			\centering
			(b) CV-WENO-SC\\[0.2cm]
			\includegraphics[width=0.95\textwidth,trim=0 0 0 0,clip]{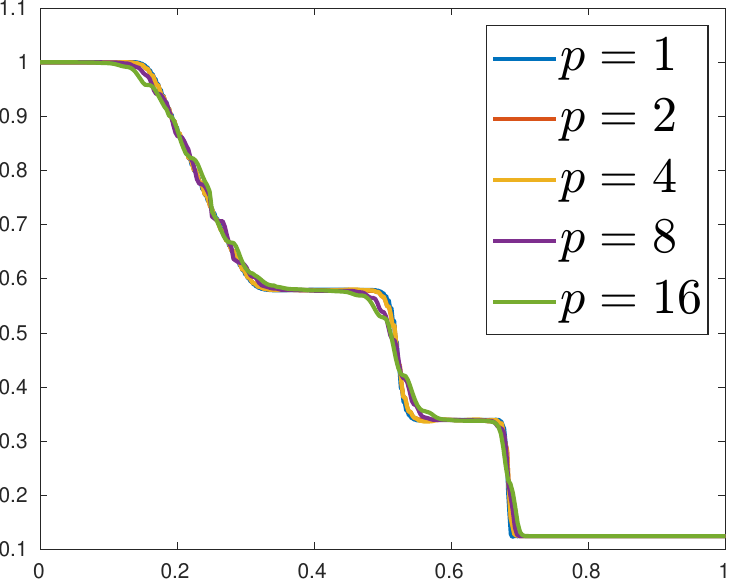}
		\end{minipage}%
		\caption{Modified Sod shock tube, density profiles $\varrho_h$ at $t=0.2$ obtained using discontinuous finite elements with $N_h=256$ DoFs per variable.}
		\label{fig:msod2}
	\end{figure}
\end{subsubsection}

\begin{subsubsection}{Woodward-Colella blast wave problem}
	The next numerical test we consider is the Woodward-Colella blast wave problem \cite{woodward1984}, a challenging benchmark for assessing the robustness of numerical schemes due to the presence of strong discontinuities in the initial pressure distribution
	\begin{align*}
		p_0(x)=\begin{cases}
			1000 & \text{if } 0<x<0.1,\\
			0.01 & \text{if } 0.1\le x\le 0.9,\\
			100 & \text{if } 0.9<x<1.
		\end{cases}
	\end{align*}
	The computational domain is $\Omega=(0,1)$ and bounded by reflective walls. The initial density is constant throughout the domain, $\varrho \equiv 1$, and the fluid is initially at rest, i.e., $v_0 \equiv 0$. Simulations are carried out up to the final time $t=0.038$, by which the two initial Riemann problems have evolved into a complex flow featuring multiple interactions between strong shock waves, rarefaction waves, and contact discontinuities. Consequently, no exact solution is available for this problem.
	
	We investigate the differences between the CV-WENO and CV-WENO-SC methods in terms of solution accuracy and shock-capturing capabilities using continuous finite elements. Figure~\ref{fig:wc} shows the numerical solutions obtained for various polynomial degrees $p$, using $N_h=1025$ degrees of freedom per variable. The CV-WENO scheme is sufficiently robust to maintain positivity of density and pressure, so that no additional flux and/or slope limiters are required to ensure invariant domain preservation (IDP). It accurately captures the right density peak, which is slightly better resolved with $p=2$ than with $p=1$, consistent with the findings in \cite{kuzmin2025a}. Notably, the shock-capturing WENO quadrature has little effect for $p\le 2$, but for $p=4$ and $p=8$ the right density peak is captured more accurately when the shock-capturing WENO quadrature is activated.
	
		\begin{figure}[t!]
		\centering
		\small
		\begin{minipage}[t]{0.48\textwidth}
			\centering
			(a) CV-WENO\\[0.2cm] 
			\includegraphics[width=0.95\textwidth,trim=0 0 0 0,clip]{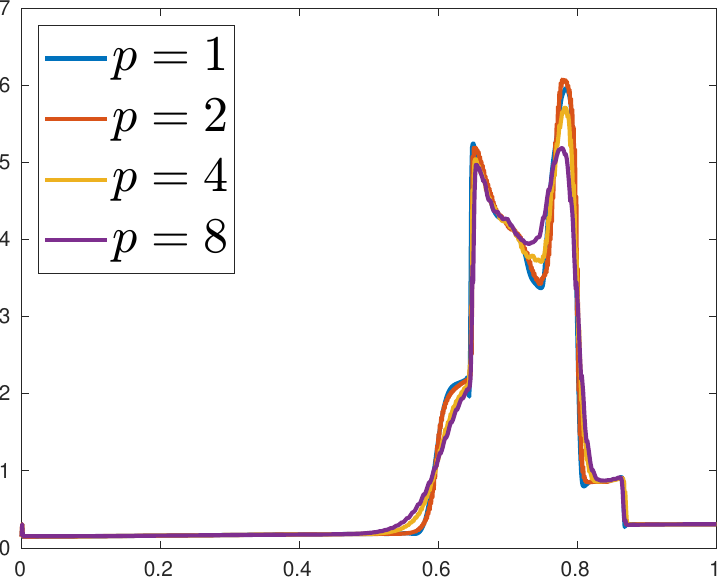}
		\end{minipage}%
		\hfill
		\begin{minipage}[t]{0.48\textwidth}
			\centering
			(b) CV-WENO-SC\\[0.2cm]
			\includegraphics[width=0.95\textwidth,trim=0 0 0 0,clip]{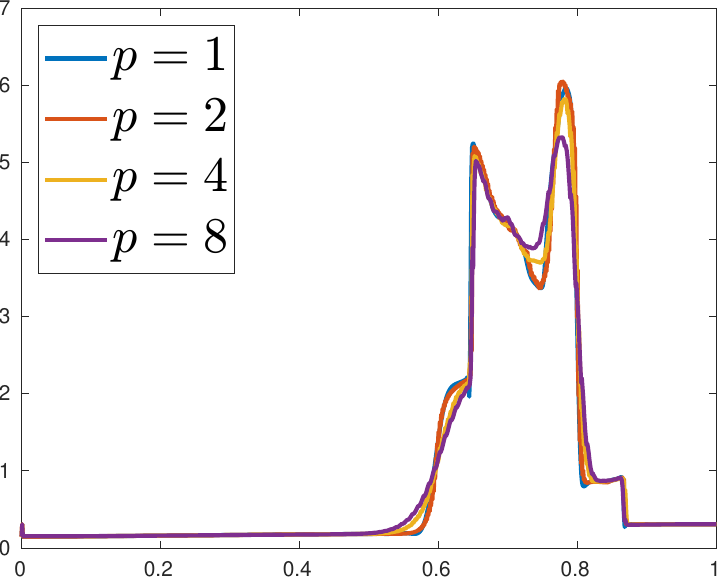}
		\end{minipage}%
		\caption{Woodward-Colella blast wave problem, density profiles $\varrho_h$ at $t=0.038$ obtained using continuous finite elements with $N_h=1025$ DoFs per variable.}
		\label{fig:wc}
	\end{figure}
\end{subsubsection}
\begin{subsubsection}{Double Mach reflection}
	In our final numerical experiment, we consider the double Mach reflection problem introduced by Woodward and Colella \cite{woodward1984}, which represents a challenging benchmark for multidimensional shock-capturing schemes. The computational domain is $\Omega=(0,4)\times(0,1)$. The pre-shock and post-shock states are given by 
	\begin{align*}
	(\varrho_L,v_{x,L},v_{y,L},p_L)&=(8.0,8.25\cos(30^\circ),-8.25\sin(30^\circ),116.5), \\ (\varrho_R,v_{x,R},v_{y,R},p_R)&=(1.4,0.0,0.0,1.0),
	\end{align*}
	 and are used to prescribe both the initial condition and the inflow boundary data. At the initial time, the post-shock state (subscript L) is imposed in the domain $\Omega_L=\{(x,y)\;|\;x<\frac{1}{6}+\frac{y}{\sqrt{3}}\}$, while the pre-shock state (subscript R) is prescribed in the complementary domain $\Omega_R = \Omega \setminus \Omega_L$. The reflecting wall corresponds to the boundary segment $y=0$ for $1/6 \le x \le 4$. No boundary conditions are required at the outflow boundary $x=4$. Along the remaining parts of the boundary, post-shock conditions are imposed for $x<\frac{1}{6}+\frac{1+20t}{\sqrt{3}}$, while pre-shock conditions are prescribed elsewhere. This construction ensures that the motion of the initial Mach $10$ shock along the top boundary is represented exactly.
	 
	 The resulting flow configuration features a Mach $10$ shock in air impinging on a reflecting wall at an initial angle of $60^\circ$. A robust numerical scheme must accurately resolve the triple-point region, while simultaneously capturing strong shocks without generating spurious oscillations.
	 
	 Figure~\ref{fig:dmr} shows snapshots of the density field at the final time $t=0.2$, computed using discontinuous finite elements on a mesh consisting of $E_h=384\cdot96$ elements with $p=2$. To prevent the occurence of negative pressures during the simulation, we additionally employ the Zhang-Shu limiter \cite{zhang2010b}. While alternative parameter choices in the WENO reconstruction may lead to IDP solutions, the scheme is not guaranteed to be IDP in general. Both the CV-WENO and CV-WENO-SC variants successfully resolve the triple-point region while introducing sufficient numerical dissipation to suppress non-physical oscillations.

	 \begin{remark}
	 	The recent work \cite{kuzmin2025a} proposes a continuous Galerkin extension of the Zhang-Shu limiting framework \cite{zhang2010b, zhang2011,zhang2012}, which can be employed to render the CG formulation of the CV-WENO scheme provably IDP. Results along these lines will be reported elsewhere.
	 \end{remark}
	 
	 \begin{figure}[t!]
	 	\centering
	 	\small
	 	\begin{minipage}[t]{0.96\textwidth}
	 		\centering
	 		(a) CV-WENO, $\varrho_h\in[1.276,22.195]$\\[0.2cm] 
	 		\includegraphics[width=0.95\textwidth,trim=0 0 0 0,clip]{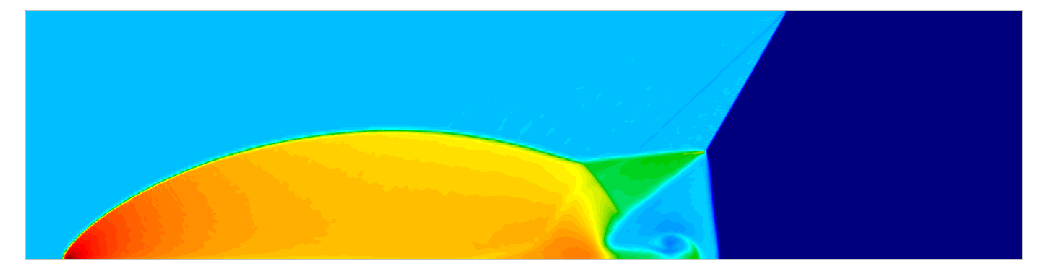}
	 	\end{minipage}%
	 	\vskip0.25cm
	 	\begin{minipage}[t]{0.96\textwidth}
	 		\centering
	 		(b) CV-WENO-SC, $\varrho_h\in [1.274, 22.423]$\\[0.2cm]
	 		\includegraphics[width=0.95\textwidth,trim=0 0 0 0,clip]{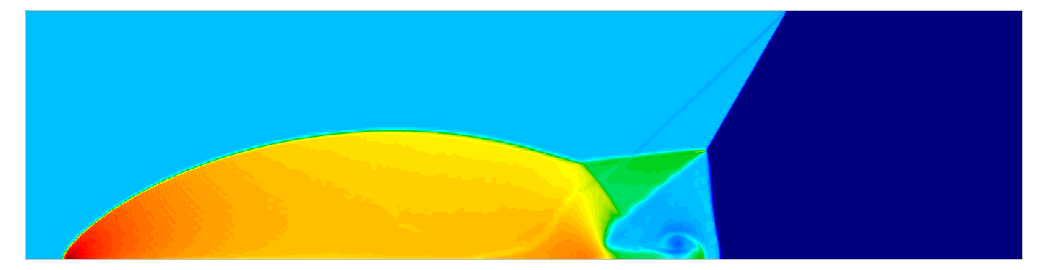}
	 	\end{minipage}%
	 	\caption{Double Mach reflection, density profiles $\varrho_h$ at $t=0.2$ obtained using discontinuous finite elements with $E_h=384\cdot96$ and $p=2$.}
	 	\label{fig:dmr}
	 \end{figure} 
\end{subsubsection}

% ---------------------------------------------------------------------------
\section{Conclusions}
\label{sec:concl}
% ---------------------------------------------------------------------------

We presented a vertex-centered WENO approach for dissipative stabilization of high-order finite element schemes, enhanced by an adaptive WENO-based scaling of quadrature weights. In our scheme, candidate polynomials of cell-level reconstructions are associated with mesh vertices rather than edges or faces, allowing information from vertex neighbors to be incorporated without increasing the number of candidates. The quadrature-weight scaling provides a computationally efficient alternative to subcell-based strategies by locally adjusting dissipation at quadrature nodes. More dissipation is applied near sharp features using a pointwise measure of solution smoothness, similarly to classical WENO schemes.
%The smoothness sensors and adaptive weights are computed on-the-fly.
The new features improve the shock-capturing capabilities of the dissipation-based WENO stabilization that we proposed in \cite{kuzmin2023a}, using its smoothness sensor to control the total amount of dissipation inside each cell. If necessary, preservation of global bounds can be enforced using the convex limiting methodology presented in \cite{kuzmin2025a}.

\section*{Acknowledgments}

The work of Dmitri Kuzmin was supported by the German Research Foundation (DFG) within the framework of the priority research program SPP 2410 under grant KU 1530/30-1. 

\bibliographystyle{abbrv}
\bibliography{references}
	
\end{document}

%% file: shared-hdr.tex
\usepackage[utf8]{inputenc}
\usepackage[ngerman,english]{babel}
\usepackage[T1]{fontenc}
\usepackage{marvosym} % Additional symbols like €

\usepackage{graphicx}
\usepackage{caption}
\usepackage{subcaption}
\usepackage{xcolor}
\usepackage{tabularx}

\setkeys{Gin}{draft=false}

\usepackage[final]{animate}

\usepackage{mathtools}
\usepackage{amssymb}
\usepackage{amsfonts}
\usepackage{amsbsy} % bold characters
\usepackage{amsopn} % allows use of \DeclareMathOperator
\usepackage{amstext} % allows use of text in math environments
\usepackage{amsxtra} % rarely used mathematical symbols
\usepackage{bbm} % for number sets, and so on
\usepackage{cancel} % cancelling in formulas
\usepackage{amscd} % commuting diagrams
\usepackage[all]{xy} % mathematical diagrams
\usepackage{esint} % more integral-types
\usepackage{empheq} % box around formulas
\usepackage{mdframed}
\usepackage{eso-pic,picture}
\usepackage{lscape}
\usepackage{blindtext}
\usepackage{listings}
\usepackage{url}
\usepackage{textcomp} % various additional text symbols
\usepackage{xspace} % fixes problems with missing spaces
\usepackage{scalerel} % use of scaling, required e.g. for reallywidehat
\usepackage{stackengine} % use of stack, required e.g. for reallywidehat
\usepackage{hhline}
\usepackage{lineno}
\usepackage{fancyhdr}
\usepackage{longtable}
\usepackage{wallpaper}
\usepackage{pdfpages}
\usepackage{siunitx} % SI units
\usepackage{hyphenat} % allows hyphenation exceptions
\usepackage{testhyphens} % provides \hyphenation

\input ushyphex
\begin{hyphenrules}{ngerman}
\end{hyphenrules}

\PassOptionsToPackage{svgnames,dvipsnames}{xcolor}
\usepackage{tcolorbox}
\tcbuselibrary{skins,breakable}

\mathcode`@="8000
{\catcode`\@=\active\gdef@{\mkern1mu}}

\newcolumntype{L}{>{$}l<{$}}
\newcolumntype{C}{>{$}c<{$}}
\newcolumntype{R}{>{$}r<{$}}

%% file: article-hdr.tex
\usepackage{lmodern}
\usepackage{microtype}

\usepackage{mathrsfs} % some special characters

\usepackage[numbers]{natbib}

\usepackage[standard,amsmath,thmmarks]{ntheorem}

\theoremstyle{break}
\theorembodyfont{\normalfont}
\theoremsymbol{{\scriptsize\ensuremath{\diamondsuit}}}

\theoremstyle{nonumberplain}
\theoremseparator{:}
\theoremsymbol{\ensuremath{\square}}
\renewtheorem{proof}{Proof}

\addto\captionsenglish{}

\usepackage[mode=multiuser]{fixme} % Provide note, warning, error, fatal.
\fxusetheme{color}

\RequirePackage[l2tabu]{nag}

\usepackage{varioref} % references with specification, where to find the reference in the text. Provided commands: \vref, \vpageref
\usepackage[colorlinks=true,citecolor=blue,
urlcolor=blue,linkcolor=blue]{hyperref} % enable clicking on a reference, this has to be included at the end
\usepackage[capitalise]{cleveref} % distinction between different types of references

\let\ccref\cref
\renewcommand{\cref}[1]{\mbox{\ccref{#1}}}
\crefname{appendix}{}{}
\crefname{listing}{Algorithm}{Algorithms}
\crefname{equation}{}{}
\crefname{table}{Table~}{Tables~}

%% file: commands.tex
\newcommand\restr[2]{{% we make the whole thing an ordinary symbol
\left.\kern-\nulldelimiterspace % automatically resize the bar with \right
#1 % the function
\vphantom{|} % pretend it's a little taller at normal size
\right|_{#2} % this is the delimiter
}}

\renewcommand{\phi}{\ensuremath{\varphi}}

\renewcommand{\epsilon}{\ensuremath{\varepsilon}}
\DeclareMathOperator{\diam}{diam}

\renewcommand{\vec}[1]{\ensuremath{\boldsymbol{#1}}}

\definecolor{green}{rgb}{0, 0.5, 0.5}
\definecolor{gray}{rgb}{0.75, 0.75, 0.75} % silver
\definecolor{slides1}{rgb}{1,0,0} % red
\definecolor{slides2}{rgb}{0,0,1} % blue
\definecolor{slides3}{rgb}{0.8, 0.5, 0.2} % bronze
\definecolor{slides4}{rgb}{0.06, 0.75, 0.99} % brightblue

\definecolor{myred}{RGB}{182,59,87}
\definecolor{myblue}{RGB}{116,135,197}
\definecolor{mygreen}{RGB}{52,152,126}
\definecolor{myorange}{RGB}{227,149,87}

\definecolor{green}{rgb}{0.0,0.6,0.0}
\definecolor{yellow}{rgb}{1,0.73,0}

\definecolor{brightgray}{RGB}{198,198,198}
\definecolor{darkgray}{RGB}{70,66,66}
\definecolor{redbrown}{RGB}{164,93,72}

\definecolor{MoSTBrightBlue}{RGB}{0,129,200}
\definecolor{MoSTDarkBlue}{RGB}{0,78,146}
\definecolor{MoSTVeryDarkBlue}{RGB}{0,58,107}

\allowdisplaybreaks % breaking of equations accross pages

\newcommand*\reallywidehat[1]{%
\savestack{\tmpbox}{\stretchto{%
\scaleto{%
\scalerel*[\widthof{\ensuremath{#1}}]{\kern-.6pt\bigwedge\kern-.6pt}%
{\rule[-\textheight/2]{1ex}{\textheight}}%WIDTH-LIMITED BIG WEDGE
}{\textheight}% 
}{0.5ex}}%
\stackon[1pt]{\ensuremath{#1}}{\tmpbox}%
}





\definecolor{keywordcolor}{rgb}{0,0.25,0.5}
\definecolor{commentcolor}{rgb}{0.2,0.5,0.2}
\definecolor{stringcolor}{rgb}{0.5,0.5,0.2}
\addto{\captionsenglish}{}

\ifdefined\vv
\renewcommand{\vv}{\ensuremath{\vec v}}
\else
\newcommand{\vv}{\ensuremath{\vec v}}
\fi

%% file: references.bib
@string{AppliedNumericalMathematics = "Appl. Numer. Math."}

@string{CommunicationsInComputationalPhysics = "Commun. Comput. Phys."}

@string{CommunicationsOnAppliedMathematicsAndComputation = "Commun. Appl. Math. Comput."}

@string{ComputerMethodsInAppliedMechanicsAndEngineering = "Comput. Methods Appl. Mech. Eng."}

@string{ComputersAndFluids = "Comput. Fluids"}

@string{ComputersAndMathematicsWithApplications = "Comput. Math. Appl."}

@string{Computing = "Computing"}

@string{ComputingAndVisualizationInScience = "Comput. Visual Sci."}

@string{ESAIMMathematicalModellingAndNumericalAnalysis = "ESAIM Math. Model. Numer. Anal."}

@string{InternationalJournalForNumericalMethodsInFluids = "Int. J. Numer. Meth. Fluids"}

@string{InternationalJournalOfComputationalFluidDynamics = "Int. J. Comput. Fluid. Dyn."}

@string{JournalOfComputationalAndAppliedMathematics = "J. Comput. Appl. Math."}

@string{JournalOfComputationalPhysics = "J. Comput. Phys."}

@string{JournalOfScientificComputing = "J. Sci. Comput."}

@string{MathematicsOfComputation = "Math. Comput."}

@string{ProceedingsOfTheRoyalSocietyAMathematicalPhysicalAndEngineeringSciences = "Proc. R. Soc. A"}

@string{Science = "Science"}

@string{SIAMJournalOnNumericalAnalysis = "SIAM J. Numer. Anal."}

@string{SIAMJournalOnScientificComputing = "SIAM J. Sci. Comput."}

@string{SIAMReview = "SIAM Rev."}

@article{adams1996,
	title={A high-resolution hybrid compact-{ENO} scheme for shock-turbulence interaction problems},
	author={Adams, Nikolaus A. and Shariff, Karim},
	journal=JournalofComputationalPhysics,
	volume={127},
	number={1},
	pages={27--51},
	year={1996},
	publisher={Elsevier}
}

@article{badia2014,
author={Badia, S and Hierro, A},
doi={10.1137/130927206},
journal=SIAMJournalOnScientificComputing,
pages={A2673--A2697},
title={On monotonicity-preserving stabilized finite element approximations of transport problems},
volume={36},
year={2014}
}

@article{badia2017a,
author={Badia, S and Bonilla, J},
doi={10.1016/j.cma.2016.09.035},
journal=ComputerMethodsInAppliedMechanicsAndEngineering,
pages={133--158},
title={Monotonicity-preserving finite element schemes based on differentiable nonlinear stabilization},
volume={313},
year={2017}
}

@article{balsara2016,
	title={An efficient class of {WENO} schemes with adaptive order},
	author={Balsara, Dinshaw S. and Garain, Sudip and Shu, Chi-Wang},
	journal=JournalofComputationalPhysics,
	volume={326},
	pages={780--804},
	year={2016},
	publisher={Elsevier}
}

@article{balsara2020,
	title={An efficient class of {WENO} schemes with adaptive order for unstructured meshes},
	author={Balsara, Dinshaw S. and Garain, Sudip and Florinski, Vladimir and Boscheri, Walter},
	journal=JournalofComputationalPhysics,
	volume={404},
	pages={109062},
	year={2020},
	publisher={Elsevier}
}

@article{barrenechea2010,
	title={Consistent local projection stabilized finite element methods},
	author={Barrenechea, Gabriel R. and Valentin, Fr{\'e}d{\'e}ric},
	journal=SIAMJournalOnNumericalAnalysis,
	volume={48},
	number={5},
	pages={1801--1825},
	year={2010},
	publisher={SIAM}
}

@article{barrenechea2013,
author={Barrenechea, G R and John, V and Knobloch, P},
doi={10.1051/m2an/2013071},
journal=ESAIMMathematicalModellingAndNumericalAnalysis,
pages={1335--1366},
volume={47},
title={A local projection stabilization finite element method  with nonlinear crosswind diffusion for convection-diffusion-reaction equations},
year={2013}
}

@article{borges2008,
	title={An improved weighted essentially non-oscillatory scheme for hyperbolic conservation laws},
	author={Borges, Rafael and Carmona, Monique and Costa, Bruno and Don, Wai Sun},
	journal=JournalofComputationalPhysics,
	volume={227},
	number={6},
	pages={3191--3211},
	year={2008},
	publisher={Elsevier}
}

@article{braack2006,
author={Braack, M. and Burman, E.},
doi={10.1137/050631227},
journal=SIAMJournalOnNumericalAnalysis,
pages={2544--2566},
title={Local projection stabilization for the {O}seen problem and its interpretation as a variational multiscale method},
volume={43},
year={2006}
}

@article{clement1975,
	title={Approximation by finite element functions using local regularization},
	author={Cl{\'e}ment, Ph},
	journal={Revue fran{\c{c}}aise d'automatique, informatique, recherche op{\'e}rationnelle. Analyse num{\'e}rique},
	volume={9},
	number={R2},
	pages={77--84},
	year={1975},
	publisher={EDP Sciences}
}

@article{codina1993,
author={Codina, R},
doi={10.1016/0045-7825(93)90213-H},
journal=ComputerMethodsInAppliedMechanicsAndEngineering,
pages={325--342},
title={A discontinuity-capturing crosswind-dissipation for the finite element solution of the convection-diffusion equation},
volume={110},
year={1993}
}

@article{codina1997,
	title={A finite element formulation for the {S}tokes problem allowing equal velocity-pressure interpolation},
	author={Codina, Ramon and Blasco, Jordi},
	journal=ComputerMethodsinAppliedMechanicsandEngineering,
	volume={143},
	number={3-4},
	pages={373--391},
	year={1997},
	publisher={Elsevier}
}

@article{codina2002,
	title={Analysis of a stabilized finite element approximation of the transient convection-diffusion-reaction equation using orthogonal subscales},
	author={Codina, Ramon and Blasco, Jordi},
	journal=ComputingAndVisualizationInScience,
	volume={4},
	number={3},
	pages={167--174},
	year={2002},
	publisher={Springer}
}

@article{friedrich1998,
author={Friedrich, O},
doi={10.1006/jcph.1998.5988},
journal=JournalOfComputationalPhysics,
pages={194--212},
title={Weighted essentially non-oscillatory schemes for the interpolation of mean values on unstructured grids},
volume={144},
year={1998}
}

@article{fu2016,
	title={A family of high-order targeted {ENO} schemes for compressible-fluid simulations},
	author={Fu, Lin and Hu, Xiangyu Y and Adams, Nikolaus A},
	journal=JournalofComputationalPhysics,
	volume={305},
	pages={333--359},
	year={2016},
	publisher={Elsevier}
}

@article{fu2017a,
	title={Targeted {ENO} schemes with tailored resolution property for hyperbolic conservation laws},
	author={Fu, Lin and Hu, Xiangyu Y and Adams, Nikolaus A},
	journal=JournalofComputationalPhysics,
	volume={349},
	pages={97--121},
	year={2017},
	publisher={Elsevier}
}

@article{glaubitz2019,
	title={Smooth and compactly supported viscous sub-cell shock capturing for discontinuous {G}alerkin methods},
	author={Glaubitz, Jan and Nogueira Jr., Alberto Costa and Almeida, Jo{\~a}o LS and Cant{\~a}o, RF and Silva, CAC},
	journal=JournalofScientificComputing,
	volume={79},
	number={1},
	pages={249--272},
	year={2019},
	publisher={Springer}
}

@misc{glvis,
	title = {{GLVis}: OpenGL Finite Element Visualization Tool},
	howpublished = {\url{glvis.org}},
	doi = {10.11578/dc.20171025.1249}
}

@article{gottlieb2001,
author={Gottlieb, S and Shu, C-W and Tadmor, E},
doi={10.1137/S003614450036757X},
journal=SIAMReview,
pages={89--112},
title={Strong stability-preserving high-order time discretization methods},
volume={43},
year={2001}
}

@article{guermond2016a,
author={Guermond, J-L and Popov, B},
doi={10.1016/j.jcp.2016.05.054},
journal=JournalOfComputationalPhysics,
pages={908--926},
title={Fast estimation from above of the maximum wave speed in the {R}iemann problem for the {E}uler equations},
volume={321},
year={2016}
}

@article{guermond2017,
author={Guermond, J-L and Popov, B},
doi={10.1137/16M1106560},
journal=SIAMJournalOnNumericalAnalysis,
pages={3120--3146},
title={Invariant domains and second-order continuous finite element approximation for scalar conservation equations},
volume={55},
year={2017}
}

@article{hajduk2020,
author={Hajduk, H and Kuzmin, D and Kolev, T and Abgrall, R},
doi={10.1016/j.cma.2019.112658},
journal=ComputerMethodsInAppliedMechanicsAndEngineering,
pages={112658},
title={Matrix-free subcell residual distribution for {B}ernstein finite element discretizations of linear advection equations},
volume={359},
year={2020}
}

@article{hajduk2020b,
author={Hajduk, H and Kuzmin, D and Kolev, T and Tomov, V and Tomas, I and Shadid, J N},
doi={10.1016/j.compfluid.2020.104451},
journal=ComputersAndFluids,
pages={104451},
title={Matrix-free subcell residual distribution for {B}ernstein finite elements: {M}onolithic limiting},
volume={200},
year={2020}
}

@article{harten1987,
author={Harten, A and Osher, S},
doi={10.1007/978-3-642-60543-7_11},
journal=SIAMJournalOnNumericalAnalysis,
pages={279--309},
title={Uniformly high-order accurate nonoscillatory schemes. {I}},
volume={24},
year={1987}
}

@article{hennemann2021,
author={Hennemann, S and Rueda-Ram{\'i}rez, A M and Hindenlang, F J and Gassner, G J},
doi={10.1016/j.jcp.2020.109935},
journal=JournalOfComputationalPhysics,
pages={109935},
title={A provably entropy stable subcell shock capturing approach for high order split form {DG} for the compressible {E}uler equations},
volume={426},
year={2021}
}

@article{henrick2005,
	title={Mapped weighted essentially non-oscillatory schemes: {A}chieving optimal order near critical points},
	author={Henrick, Andrew K and Aslam, Tariq D and Powers, Joseph M},
	journal=JournalofComputationalPhysics,
	volume={207},
	number={2},
	pages={542--567},
	year={2005},
	publisher={Elsevier}
}

@article{hill2004,
	title={Hybrid tuned center-difference-{WENO} method for large eddy simulations in the presence of strong shocks},
	author={Hill, David J and Pullin, Dale I},
	journal=JournalofComputationalPhysics,
	volume={194},
	number={2},
	pages={435--450},
	year={2004},
	publisher={Elsevier}
}

@article{hu2010,
	title={An adaptive central-upwind weighted essentially non-oscillatory scheme},
	author={Hu, XY and Wang, Q and Adams, Nikolaus A},
	journal=JournalofComputationalPhysics,
	volume={229},
	number={23},
	pages={8952--8965},
	year={2010},
	publisher={Elsevier}
}

@article{hughes1986,
author={Hughes, T J R and Mallet, M},
doi={10.1016/0045-7825(86)90153-2},
journal=ComputerMethodsInAppliedMechanicsAndEngineering,
pages={329--336},
title={A new finite element formulation for computational fluid dynamics: {IV}. {A} discontinuity-capturing operator for multidimensional advective-diffusive systems},
volume={58},
year={1986}
}

@article{hughes1987,
  title={Recent progress in the development and understanding of {SUPG} methods with special reference to the compressible {E}uler and {N}avier-{S}tokes equations},
  author={Hughes, Thomas J R},
  journal=InternationalJournalForNumericalMethodsInFluids,
  volume={7},
  number={11},
  pages={1261--1275},
  year={1987}
}

@inproceedings{jameson1981,
author={Jameson, A and Schmidt, W and Turkel, E},
booktitle={14th Fluid and Plasma Dynamics Conference},
doi={10.2514/6.1981-1259},
publisher={American Institute of Aeronautics and Astronautics},
title={Numerical solution of the {E}uler equations by finite volume methods using {R}unge {K}utta time stepping schemes},
year={1981}
}

@article{jiang1996,
author={Jiang, G-S and Shu, C-W},
doi={10.1006/jcph.1996.0130},
journal=JournalOfComputationalPhysics,
pages={202--228},
title={Efficient implementation of weighted {ENO} schemes},
volume={126},
year={1996}
}

@article{john2006,
author={John, V and Kaya, S and Layton, W},
doi={10.1016/j.cma.2005.10.006},
journal=ComputerMethodsInAppliedMechanicsAndEngineering,
pages={4594--4603},
title={A two-level variational multiscale method for convection-dominated convection-diffusion equations},
volume={195},
year={2006}
}

@article{johnson1990,
author={Johnson, C and Szepessy, A and Hansbo, P},
doi={10.2307/2008684},
journal=MathematicsOfComputation,
volume={54},
pages={107--129},
title={On the convergence of shock-capturing streamline diffusion finite element methods for hyperbolic conservation laws},
year={1990}
}

@article{knobloch2009,
	title={Local projection stabilization for advection--diffusion--reaction problems: One-level vs. two-level approach},
	author={Knobloch, Petr and Lube, Gert},
	journal=AppliedNumericalMathematics,
	volume={59},
	number={12},
	pages={2891--2907},
	year={2009},
	publisher={Elsevier}
}

@article{knobloch2010a,
	title={A generalization of the local projection stabilization for convection-diffusion-reaction equations},
	author={Knobloch, Petr},
	journal=SIAMJournalonNumericalAnalysis,
	volume={48},
	number={2},
	pages={659--680},
	year={2010},
	publisher={SIAM}
}

@article{kurganov2007,
author={Kurganov, A and Petrova, G and Popov, B},
doi={10.1137/040614189},
journal=SIAMJournalOnScientificComputing,
pages={2381--2401},
title={Adaptive semidiscrete central-upwind schemes for nonconvex hyperbolic conservation laws},
volume={29},
year={2007}
}

@article{kuzmin2010,
author={Kuzmin, D},
doi={10.1016/j.cam.2009.05.028},
journal=JournalOfComputationalAndAppliedMathematics,
pages={3077--3085},
title={A vertex-based hierarchical slope limiter for p-adaptive discontinuous {G}alerkin methods},
volume={233},
year={2010}
}

@article{kuzmin2020,
author={Kuzmin, D},
doi={10.1016/j.cma.2019.112804},
journal=ComputerMethodsInAppliedMechanicsAndEngineering,
pages={112804},
title={Monolithic convex limiting for continuous finite element discretizations of hyperbolic conservation laws},
volume={361},
year={2020}
}

@article{kuzmin2020a,
author={Kuzmin, D and Quezada de Luna, M},
doi={10.1016/j.jcp.2020.109411},
journal=JournalOfComputationalPhysics,
pages={109411},
title={Subcell flux limiting for high-order {B}ernstein finite element discretizations of scalar hyperbolic conservation laws},
volume={411},
year={2020}
}

@article{kuzmin2020f,
author={Kuzmin, D and Quezada de Luna, M},
doi={10.1016/j.compfluid.2020.104742},
journal=ComputersAndFluids,
pages={104742},
title={Entropy conservation property and entropy stabilization of high-order continuous {G}alerkin approximations to scalar conservation laws},
volume={213},
year={2020}
}

@article{kuzmin2023a,
author={Kuzmin, D and Vedral, J},
doi={10.1016/j.jcp.2023.112153},
journal=JournalOfComputationalPhysics,
pages={112153},
title={Dissipation-based {WENO} stabilization of high-order finite element methods for scalar conservation laws},
volume={487},
year={2023}
}

@article{kuzmin2025a,
	title={A matrix-free convex limiting framework for continuous {G}alerkin methods with nonlinear stabilization},
	author={Kuzmin, Dmitri and Hajduk, Hennes and Vedral, Joshua},
	journal={arXiv preprint arXiv:2509.04673},
	year={2025}
}

@article{leveque1996,
author={LeVeque, R J},
doi={10.1137/0733033},
journal=SIAMJournalOnNumericalAnalysis,
pages={627--665},
title={High-resolution conservative algorithms for advection in incompressible flow},
volume={33},
year={1996}
}

@article{levy1999,
	title={Central {WENO} schemes for hyperbolic systems of conservation laws},
	author={Levy, Doron and Puppo, Gabriella and Russo, Giovanni},
	journal=ESAIMMathematicalModellingAndNumericalAnalysis,
	volume={33},
	number={3},
	pages={547--571},
	year={1999},
	publisher={EDP Sciences}
}

@article{levy2000,
	title={Compact central {WENO} schemes for multidimensional conservation laws},
	author={Levy, Doron and Puppo, Gabriella and Russo, Giovanni},
	journal=SIAMJournalonScientificComputing,
	volume={22},
	number={2},
	pages={656--672},
	year={2000},
	publisher={SIAM}
}

@article{li2020,
	title={An efficient low-dissipation hybrid central/{WENO} scheme for compressible flows},
	author={Li, Liang and Wang, Hong-Bo and Zhao, Guo-Yan and Sun, Ming-Bo and Xiong, Da-Peng and Tang, Tao},
	journal=InternationalJournalofComputationalFluidDynamics,
	volume={34},
	number={10},
	pages={705--730},
	year={2020},
	publisher={Taylor \& Francis}
}

@article{liu1994,
author={Liu, X-D and Osher, S and Chan, T},
doi={10.1006/jcph.1994.1187},
journal=JournalOfComputationalPhysics,
pages={200--212},
title={Weighted essentially non-oscillatory schemes},
volume={115},
year={1994}
}

@article{lohmann2017,
author={Lohmann, C and Kuzmin, D and Shadid, J N and Mabuza, S},
doi={10.1016/j.jcp.2017.04.059},
journal=JournalOfComputationalPhysics,
pages={151--186},
title={Flux-corrected transport algorithms for continuous {G}alerkin methods based on high order {B}ernstein finite elements},
volume={344},
year={2017}
}

@article{luo2007,
author={Luo, H and Baum, J D and L{\"o}hner, R},
doi={10.1016/j.jcp.2006.12.017},
journal=JournalOfComputationalPhysics,
pages={686--713},
title={A {H}ermite {WENO}-based limiter for discontinuous {G}alerkin method on unstructured grids},
volume={225},
year={2007}
}

@article{luo2008,
author={Luo, H and Baum, J D and L{\"o}hner, R},
doi={10.1016/j.jcp.2008.06.035},
journal=JournalOfComputationalPhysics,
pages={8875--8893},
title={A discontinuous {G}alerkin method based on a {T}aylor basis for the compressible flows on arbitrary grids},
volume={227},
year={2008}
}

@article{markert2023,
	title={A sub-element adaptive shock capturing approach for discontinuous {G}alerkin methods},
	author={Markert, Johannes and Gassner, Gregor and Walch, Stefanie},
	journal=CommunicationsonAppliedMathematicsandComputation,
	volume={5},
	number={2},
	pages={679--721},
	year={2023},
	publisher={Springer}
}

@misc{mfem,
	key = {mfem},
	title = {{MFEM}: Modular Finite Element Methods {[Software]}},
	howpublished = {\url{mfem.org}},
	doi = {10.11578/dc.20171025.1248}
}

@article{mfem2021,
	title = {{MFEM}: A Modular Finite Element Methods Library},
	author = {R. Anderson and J. Andrej and A. Barker and J. Bramwell and J.-S. Camier and
	J. Cerveny and V. Dobrev and Y. Dudouit and A. Fisher and Tz. Kolev and W. Pazner and
	M. Stowell and V. Tomov and I. Akkerman and J. Dahm and D. Medina and S. Zampini},
	journal = {Computers \& Mathematics with Applications},
	doi = {10.1016/j.camwa.2020.06.009},
	volume = {81},
	pages = {42-74},
	year = {2021}
}

@article{mfem2024,
	title = {High-Performance Finite Elements with {MFEM}},
	author = {J. Andrej and N. Atallah and J.-P. Bäcker and  J.-S. Camier and
	D. Copeland and V. Dobrev and Y. Dudouit and T. Duswald and
	B. Keith and D. Kim and T. Kolev and B. Lazarov and  K. Mittal and
	W. Pazner and S. Petrides and S. Shiraiwa and M. Stowell and V. Tomov},
	journal={The International Journal of High Performance Computing Applications},
	volume={38},
	number={5},
	pages={447-467},
	year={2024},
	publisher={SAGE Publications Sage UK: London, England}
}

@article{navas-montilla2024,
	title={Exploring the potential of {TENO} and {WENO} schemes for simulating under-resolved turbulent flows in the atmosphere using {E}uler equations},
	author={Navas-Montilla, Adri{\'a}n and Guallart, Javier and Sol{\'a}n-Fustero, Pablo and Garc{\'\i}a-Navarro, Pilar},
	journal=ComputersAndFluids,
	volume={280},
	pages={106349},
	year={2024},
	publisher={Elsevier}
}

@article{nazarov2013,
author={Nazarov, M},
doi={10.1016/j.camwa.2012.11.003},
journal=ComputersAndMathematicsWithApplications,
pages={616--626},
title={Convergence of a residual based artificial viscosity finite element method},
volume={65},
year={2013}
}

@incollection{persson2006,
author={Persson, P-O and Peraire, J},
booktitle={44th AIAA Aerospace Sciences Meeting and Exhibit},
doi={10.2514/6.2006-112},
publisher={American Institute of Aeronautics and Astronautics},
title={Sub-cell shock capturing for discontinuous {G}alerkin methods},
year={2006}
}

@article{pirozzoli2002,
	title={Conservative hybrid compact-{WENO} schemes for shock-turbulence interaction},
	author={Pirozzoli, Sergio},
	journal=JournalofComputationalPhysics,
	volume={178},
	number={1},
	pages={81--117},
	year={2002},
	publisher={Elsevier}
}

@article{qiu2002,
	title={On the construction, comparison, and local characteristic decomposition for high-order central {WENO} schemes},
	author={Qiu, Jianxian and Shu, Chi-Wang},
	journal=JournalofComputationalPhysics,
	volume={183},
	number={1},
	pages={187--209},
	year={2002},
	publisher={Elsevier}
}

@article{qiu2004,
	title={Hermite {WENO} schemes and their application as limiters for {R}unge--{K}utta discontinuous {G}alerkin method: {O}ne-dimensional case},
	author={Qiu, Jianxian and Shu, Chi-Wang},
	journal=JournalofComputationalPhysics,
	volume={193},
	number={1},
	pages={115--135},
	year={2004},
	publisher={Elsevier}
}

@article{qiu2005,
author={Qiu, J and Shu, C-W},
doi={10.1137/S1064827503425298},
journal=SIAMJournalOnScientificComputing,
pages={907--929},
title={Runge--{K}utta discontinuous {G}alerkin method using {WENO} limiters},
volume={26},
year={2005}
}

@article{qiu2005b,
	title={Hermite {WENO} schemes and their application as limiters for {R}unge--{K}utta discontinuous {G}alerkin method {II}: {T}wo dimensional case},
	author={Qiu, Jianxian and Shu, Chi-Wang},
	journal=ComputersAndFluids,
	volume={34},
	number={6},
	pages={642--663},
	year={2005},
	publisher={Elsevier}
}

@book{quarteroni1994,
author={Quarteroni, A and Valli, A},
doi={10.1007/978-3-540-85268-1},
publisher={Springer},
title={Numerical approximation of partial differential equations},
year={1994}
}

@article{rueda-ramirez2022,
author={Rueda-Ram{\'i}rez, A M and Pazner, W and Gassner, G J},
doi={10.1016/j.compfluid.2022.105627},
journal=ComputersAndFluids,
pages={105627},
title={Subcell limiting strategies for discontinuous {G}alerkin spectral element methods},
volume={247},
year={2022}
}

@article{scott1990,
	title={Finite element interpolation of nonsmooth functions satisfying boundary conditions},
	author={Scott, L Ridgway and Zhang, Shangyou},
	journal=MathematicsofComputation,
	volume={54},
	number={190},
	pages={483--493},
	year={1990}
}

@article{shu1988,
author={Shu, C-W and Osher, S},
doi={10.1016/0021-9991(88)90177-5},
journal=JournalOfComputationalPhysics,
pages={439--471},
title={Efficient implementation of essentially non-oscillatory shock-capturing schemes},
volume={77},
year={1988}
}

@incollection{shu1989,
author={Shu, C-W and Osher, S},
booktitle={Upwind and High-Resolution Schemes},
doi={10.1007/978-3-642-60543-7_14},
pages={328--374},
publisher={Springer},
title={Efficient implementation of essentially non-oscillatory shock-capturing schemes,
{II}},
year={1989}
}

@article{sod1978,
author={Sod, G A},
doi={10.1016/0021-9991(78)90023-2},
journal=JournalOfComputationalPhysics,
pages={1--31},
title={A survey of several finite difference methods for systems of nonlinear hyperbolic conservation laws},
volume={27},
year={1978}
}

@inproceedings{sonntag2014,
author={Sonntag, M and Munz, C-D},
booktitle={Finite Volumes for Complex Applications VII - Elliptic, Parabolic and Hyperbolic Problems},
doi={10.1007/978-3-319-05591-6_96},
pages={945--953},
publisher={Springer},
title={Shock capturing for discontinuous {G}alerkin methods using finite volume subcells},
year={2014}
}

@book{toro2013,
	title={Riemann solvers and numerical methods for fluid dynamics: {a} practical introduction},
	author={Toro, Eleuterio F.},
	year={2013},
	publisher={Springer Science \& Business Media}
}

@article{tsoutsanis2019,
	title={Stencil selection algorithms for {WENO} schemes on unstructured meshes},
	author={Tsoutsanis, Panagiotis},
	journal={J. Comput. Phys.: X},
	volume={4},
	pages={100037},
	year={2019},
	publisher={Elsevier}
}

@article{vedral2025,
author={Vedral, J and Rupp, A and Kuzmin, D},
doi={10.1016/j.apnum.2024.12.008},
journal=AppliedNumericalMathematics,
pages={64--81},
title={Strongly consistent low-dissipation {WENO} schemes for finite elements},
volume={210},
year={2025}
}

@article{vedral-arxiv,
author={Vedral, J.},
journal={Preprint arXiv:2309.12019},
title={Dissipative {WENO} stabilization of high-order discontinuous {G}alerkin methods for hyperbolic problems},
year={2023}
}

@article{vilar2024,
author={Vilar, F and Abgrall, R},
doi={10.1137/22M1542696},
journal=SIAMJournalOnScientificComputing,
pages={A851--A883},
title={A posteriori local subcell correction of high-order discontinuous {G}alerkin scheme for conservation laws on two-dimensional unstructured grids},
volume={46},
year={2024}
}

@article{wan2012,
	title={Robustness of the hybrid {DRP}-{WENO} scheme for shock flow computations},
	author={Wan, Zhen-Hua and Zhou, Lin and Sun, De-Jun},
	journal=InternationalJournalForNumericalMethodsInFluids,
	volume={70},
	number={8},
	pages={985--1003},
	year={2012},
	publisher={Wiley Online Library}
}

@article{woodward1984,
author={Woodward, P and Colella, P},
doi={10.1016/0021-9991(84)90142-6},
journal=JournalOfComputationalPhysics,
pages={115--173},
title={The numerical simulation of two-dimensional fluid flow with strong shocks},
volume={54},
year={1984}
}

@article{zhang2010b,
author={Zhang, X and Shu, C-W},
doi={10.1016/j.jcp.2010.08.016},
journal=JournalOfComputationalPhysics,
pages={8918--8934},
title={On positivity-preserving high order discontinuous {G}alerkin schemes for compressible {E}uler equations on rectangular meshes},
volume={229},
year={2010}
}

@article{zhang2011,
author={Zhang, X and Shu, C-W},
doi={10.1098/rspa.2011.0153},
journal=ProceedingsOfTheRoyalSocietyAMathematicalPhysicalAndEngineeringSciences,
pages={2752--2776},
title={Maximum-principle-satisfying and positivity-preserving high-order schemes for conservation laws: {S}urvey and new developments},
volume={467},
year={2011}
}

@article{zhang2012,
author={Zhang, X and Xia, Y and Shu, C-W},
doi={10.1007/s10915-011-9472-8},
journal=JournalOfScientificComputing,
pages={29--62},
title={Maximum-principle-satisfying and positivity-preserving high order discontinuous {G}alerkin schemes for conservation laws on triangular meshes},
volume={50},
year={2012}
}

@article{zhong2013,
author={Zhong, X and Shu, C-W},
doi={10.1016/j.jcp.2012.08.028},
journal=JournalOfComputationalPhysics,
pages={397--415},
title={A simple weighted essentially nonoscillatory limiter for {R}unge--{K}utta discontinuous {G}alerkin methods},
volume={232},
year={2013}
}

@article{zhu2008,
  title={Runge--{K}utta discontinuous {G}alerkin method using {WENO} limiters {II}: {U}nstructured meshes},
  author={Zhu, Jun and Qiu, Jianxian and Shu, Chi-Wang and Dumbser, Michael},
  journal=JournalOfComputationalPhysics,
  volume={227},
  number={9},
  pages={4330--4353},
  year={2008}
}

@article{zhu2009,
author={Zhu, J and Qiu, J},
doi={10.1007/s10915-009-9271-7},
journal=JournalOfScientificComputing,
pages={293--321},
title={Hermite {WENO} schemes and their application as limiters for {R}unge--{K}utta discontinuous {G}alerkin method, {III}: {U}nstructured meshes},
volume={39},
year={2009}
}

@article{zhu2016,
	title={A new fifth order finite difference {WENO} scheme for solving hyperbolic conservation laws},
	author={Zhu, Jun and Qiu, Jianxian},
	journal=JournalofComputationalPhysics,
	volume={318},
	pages={110--121},
	year={2016},
	publisher={Elsevier}
}

@article{zhu2017,
author={Zhu, J and Zhong, X and Shu, C-W and Qiu, J},
doi={10.4208/cicp.221015.160816a},
journal=CommunicationsInComputationalPhysics,
pages={623--649},
title={Runge--{K}utta discontinuous {G}alerkin method with a simple and compact {H}ermite {WENO} limiter on unstructured meshes},
volume={21},
year={2017}
}

@article{zhu2018,
	title={New finite volume weighted essentially nonoscillatory schemes on triangular meshes},
	author={Zhu, Jun and Qiu, Jianxian},
	journal=SIAMJournalonScientificComputing,
	volume={40},
	number={2},
	pages={A903--A928},
	year={2018},
	publisher={SIAM}
}

@article{zou2025,
	title={High-Order Finite-Volume Multi-resolution {WENO} Schemes on Mixed-Element Unstructured Meshes for Simulating Two-Dimensional Compressible Flows},
	author={Zou, Zimai and Zhu, Jun and Wang, Chunwu},
	journal=CommunicationsonAppliedMathematicsandComputation,
	pages={1--22},
	year={2025},
	publisher={Springer}
}

@article{olshanskii2025,
  title={Gradient-penalty stabilization of sharp and diffuse interface formulations in unfitted {N}itsche finite element methods},
  author={Olshanskii, M and B{\"a}cker, J-P and Kuzmin, D},
  journal={Preprint arXiv:2501.16594},
  year={2025}
}
